\magnification=1200
\input amssym.def
\input amssym.tex

\hsize=13.5truecm
\baselineskip=16truept

\def\today{${\scriptscriptstyle\number\day-\number\month-\number\year}$}
\footline={{\hfil\rm\the\pageno\hfil${\scriptscriptstyle\rm\jobname}$\ \ \today}}

\font\secbf=cmb10 scaled 1200
\font\eightrm=cmr8
\font\sixrm=cmr6

\font\eighti=cmmi8

\font\sixi=cmmi6
\skewchar\eighti='177 \skewchar\sixi='177

\font\eightsy=cmsy8
\font\sixsy=cmsy6
\skewchar\eightsy='60 \skewchar\sixsy='60

\font\eightit=cmti8

\font\eightbf=cmbx8
\font\sixbf=cmbx6

\let\sc=\tensc

\font\eightsc=cmcsc10 scaled 800
\font\secbf=cmb10 scaled 1200
\font\subsecfont=cmb10 scaled \magstephalf
\font\amb=cmmib10

\font\ambi=cmmib10 scaled 700

\newfam\mbfam 

\textfont\mbfam\amb \scriptfont\mbfam\ambi


\def\aa{\def\rm{\fam0\eightrm}%
  \textfont0=\eightrm \scriptfont0=\sixrm \scriptscriptfont0=\fiverm
  \textfont1=\eighti \scriptfont1=\sixi \scriptscriptfont1=\fivei
  \textfont2=\eightsy \scriptfont2=\sixsy \scriptscriptfont2=\fivesy
  \textfont3=\tenex \scriptfont3=\tenex \scriptscriptfont3=\tenex
  \def\sc{\eightsc}
  \def\it{\fam\itfam\eightit}%
  \textfont\itfam=\eightit
  \def\bf{\fam\bffam\eightbf}%
  \textfont\bffam=\eightbf \scriptfont\bffam=\sixbf
   \scriptscriptfont\bffam=\fivebf
  \normalbaselineskip=9.7pt
  \setbox\strutbox=\hbox{\vrule height7pt depth2.6pt width0pt}%
  \normalbaselines\rm}

\def\Proof{\vskip12pt\noindent{\bf Proof.\ }}

\def\Def#1{\vskip12pt\noindent{\bf Definition #1}}
\def\Remark#1{\vskip12pt\noindent{\bf Remark #1}}

\def\m@th{\mathsurround=0pt}

\def\cc#1{\hbox to .88\hsize{$\displaystyle\hfil{#1}\hfil$}\cr}
\def\lc#1{\hbox to .88\hsize{$\displaystyle{#1}\hfill$}\cr}
\def\rc#1{\hbox to .88\hsize{$\displaystyle\hfill{#1}$}\cr}

\def\eqal#1{\null\,\vcenter{\openup\jot\m@th
  \ialign{\strut\hfil$\displaystyle{##}$&&$\displaystyle{{}##}$\hfil
      \crcr#1\crcr}}\,}

\def\section#1{\vskip 22pt plus6pt minus2pt\penalty-400
        {{\secbf
        \noindent#1\rightskip=0pt plus 1fill\par}}
        \par\vskip 12pt plus5pt minus 2pt
        \penalty 1000}

\def\subsection#1{\vskip 20pt plus6pt minus2pt\penalty-400
        {{\subsecfont
        \noindent#1\rightskip=0pt plus 1fill\par}}
        \par\vskip 8pt plus5pt minus 2pt
        \penalty 1000}

\def\subsubsection#1{\vskip 18pt plus6pt minus2pt\penalty-400
        {{\subsecfont
        \noindent#1}}
        \par\vskip 7pt plus5pt minus 2pt
        \penalty 1000}

\def\center#1{{\begingroup \leftskip=0pt plus 1fil\rightskip=\leftskip
\parfillskip=0pt \spaceskip=.3333em \xspaceskip=.5em \pretolerance 9999
\tolerance 9999 \parindent 0pt \hyphenpenalty 9999 \exhyphenpenalty 9999
\par #1\par\endgroup}}

\def\\{\hfill\break}

\def\mida#1{{{\null\kern-4.2pt\left\bracevert\vbox to 6pt{}\!\hbox{$#1$}\!\right\bracevert\!\!}}}
\def\midy#1{{{\null\kern-4.2pt\left\bracevert\!\!\hbox{$\scriptstyle{#1}$}\!\!\right\bracevert\!\!}}}

\def\diagint{{\raise1.5pt\hbox{$\scriptscriptstyle\diagup$}\hskip-8.7pt\intop}}

\def\divv{{\rm div}\,}
\def\rot{{\rm rot}\,}
\def\const{{\rm const}}

\def\krop{{\hbox{$\scriptscriptstyle\bullet$}}}

\def\krop{{\hbox{$\scriptscriptstyle\bullet$}}}

\def\D{{\Bbb D}}

\def\I{{\Bbb I}}
\def\R{{\Bbb R}}
\def\N{{\Bbb N}}
\def\T{{\Bbb T}}

\def\krop{{\hbox{$\scriptscriptstyle\bullet$}}}
\center{\secbf Long time existence of regular solutions to non-homogeneous 
Navier-Stokes equations}
\vskip1.5cm

\centerline{\bf Wojciech M. Zaj\c aczkowski\footnote{$^*$}{\aa The author is 
partially supported by Polish Grant $\rm NN\,201\,396\,937$}}

\vskip1cm
\noindent
Institute of Mathematics, Polish Academy of Sciences,\\
\'Sniadeckich 8, 00-956 Warsaw, Poland\\
E-mail:wz@impan.pl;\\
Institute of Mathematics and Cryptology, Cybernetics Faculty,\\
Military University of Technology, Kaliskiego 2,\\
00-908 Warsaw, Poland
\vskip0.8cm


\vskip1.5cm

\noindent
{\bf Abstract.} We consider the motion of incompressible viscous 
non-homo\-geneous fluid described by the Navier-Stokes equations in a bounded 
cylinder $\Omega$ under boundary slip conditions. Assume that the $x_3$-axis 
is the axis of the cylinder. Let $\varrho$ be the density of the fluid, 
$v$ -- the velocity and $f$ the external force field. Assuming that 
quantities $\nabla\varrho(0)$, $\partial_{x_3}v(0)$, $\partial_{x_3}f$, 
$f_3|_{\partial\Omega}$ are sufficiently small in some norms we prove large 
time regular solutions such that $v\in H^{2+s,1+s/2}(\Omega\times(0,T))$, 
$\nabla p\in H^{s,s/2}(\Omega\times(0,T))$, ${1\over2}<s<1$ without any 
restriction on the existence time $T$. The proof is divided into two parts. 
First an a priori estimate is shown. Next the existence follows from the 
Leray-Schauder fixed point theorem.
\vfil\eject

\section{1. Introduction}

We consider the motion of a viscous non-homogeneous incompressible fluid 
described by the following system
$$\eqal{
&\varrho(v_{,t}+v\cdot\nabla v)-\divv\T(v,p)=\varrho f\quad &{\rm in}\ \ 
\Omega^T=\Omega\times(0,T),\cr
&\divv v=0\quad &{\rm in}\ \ \Omega^T,\cr
&\varrho_{,t}+v\cdot\nabla\varrho=0\quad &{\rm in}\ \ \Omega^T,\cr
&v\cdot\bar n=0\quad &{\rm on}\ \ S^T=S\times(0,T),\cr
&\bar n\cdot\T(v,p)\cdot\bar\tau_\alpha+\delta_{i1}\gamma v\cdot
\bar\tau_\alpha=0,\ \ \alpha=1,2,\quad &{\rm on}\ \ S_i^T,\ \ i=1,2,\cr
&v|_{t=0}=v_0\quad &{\rm in}\ \ \Omega,\cr
&\varrho|_{t=0}=\varrho_0\quad &{\rm in}\ \ \Omega,\cr}
\leqno(1.1)
$$
where $\Omega\subset\R^3$ is a bounded cylindrical domain, $S=\partial\Omega$ 
is the boundary of $\Omega$, $v=(v_1(x,t),v_2(x,t),v_3(x,t))\in\R^3$ 
is the velocity of the fluid, $\varrho=\varrho(x,t)\in\R_+$ the density, 
$p=p(x,t)\in\R$ the pressure, $f=(f_1(x,t),\break
f_2(x,t),f_3(x,t))\in\R^3$ the 
external force field and $x=(x_1,x_2,x_3)$ are the Cartesian coordinates. 
Moreover, $\bar n$ is the unit outward vector normal to $S$, $\bar\tau_i$, 
$i=1,2$, are tangent vectors to $S$ and the dot denotes the scalar product in 
$\R^3$. Finally $\gamma\ge0$ is the constant slip coefficient and 
$\delta_{ij}$ is the Kronecker $\delta$.

By $\T(v,p)$ we denote the stress tensor of the form
$$
\T(v,p)=\nu\D(v)-p\I,
\leqno(1.2)
$$
where $\nu$ is the constant viscosity coefficient, $\I$ is the unit matrix 
and $\D(v)$ is the dilatation tensor of the form
$$
\D(v)=\{v_{i,x_j}+v_{j,x_i}\}_{i,j=1,2,3}.
\leqno(1.3)
$$
$\Omega$ is cylindrical domain parallel to the $x_3$-axis with arbitrary 
cross section.

\noindent
We assume that $S=S_1\cup S_2$, where $S_1$ is the part of the boundary 
parallel to the $x_3$-axis and $S_2$ is perpendicular to $x_3$. Hence,
$$\eqal{
&S_1=\{x\in\R^3:\ \psi(x_1,x_2)=c_0,\ -a<x_3<a\},\cr
&S_2=\{x\in\R^3:\ \psi(x_1,x_2)<c_0,\ x_3\in\{-a,a\}\},\cr}
\leqno(1.4)
$$
where $a$ and $c_0$ are given positive numbers and $\psi(x_1,x_2)=c_0$ 
describes a sufficiently smooth closed curve in the plane $x_3=\const$.

\noindent
Now we formulate the main result of this paper

\proclaim Theorem A. 
Assume that
\item{1.} $v_0\in H^{1+s}(\Omega)$, $v_{0,x_3}\in H^1(\Omega)$, 
$s\in(1/2,1)$;
\item{2.} $\varrho_0\in W_q^2(\Omega)$, $3<q\le{3\over3/2-s}$, 
$s\in(1/2,1)$;
\item{3.} there exist positive constants $0<\varrho_*<\varrho^*$ such that 
$\varrho_*\le\varrho_0\le\varrho^*$;
\item{4.} $f\in H^{s,s/2}(\Omega^T)$, $f_{,x_3}\in L_2(\Omega^T)$, 
$s\in(1/2,1)$;
\item{5.} the quantity
$$\eqal{
X(T)&=\|\varrho_{0,x}\|_{W_q^1(\Omega)}+
\|f_{,x_3}\|_{L_2(0,T;L_{6/5}(\Omega))}\cr
&\quad+\|v_{0,x_3}\|_{L_2(\Omega)}+\|f_3\|_{L_2(0,T;L_{4/3}(S_2))}\le\delta\cr}
\leqno(1.4)
$$
where $\delta$ is sufficiently small.

{\noindent\sl 
Then there exists a unique solution to problem (1.1) such that \\
$v\in H^{2+s,1+s/2}(\Omega^T)$, $v_{,x_3}\in H^{2,1}(\Omega^T)$, 
$\nabla p\in H^{s,s/2}(\Omega^T)$, $\nabla p_{,x_3}\in L_2(\Omega^T)$
$$\eqal{
&\|v\|_{H^{2+s,1+s/2}(\Omega^T)}+\|v_{,x_3}\|_{H^{2,1}(\Omega^T)}+
\|\nabla p\|_{H^{s,s/2}(\Omega^T)}\cr
&\quad+\|\nabla p_{,x_3}\|_{L_2(\Omega^T)}\le\varphi(\varrho_*,\varrho^*,N)\cr}
\leqno(1.5)
$$
where $\varphi$ is an increasing positive function and
$$\eqal{
N&=\|f\|_{H^{s,s/2}(\Omega^T)}+\|f_{,x_3}\|_{L_2(\Omega^T)}+
\|f\|_{L_2(0,T;W_{6/5}^1(\Omega))}\cr
&\quad+\|v_0\|_{H^{1+s}(\Omega)}+\|v_{0,x_3}\|_{H^1(\Omega)}.\cr}
$$}

\noindent
The result formulated in Theorem A describes a long time existence of 
solutions to problem (1.1) because the smallness condition (1.4) contains 
at most time integral norms of $f$. 

\noindent
The aim of this paper is to prove long time existence of regular solutions 
to problem (1.1) such that there is no restriction on the magnitudes of the 
external force, the initial velocity and the density. The aim is covered by 
the smallness restriction (1.4) because it contains derivatives of the initial 
density and derivatives with respect to $x_3$ of the initial velocity and 
the external force. This kind of restrictions suggests that our solution 
remains close to two-dimensional solutions of incompressible Navier-Stokes 
equations because the initial density is close to a constant but the initial 
velocity and the external force change a little in the $x_3$-direction. 
In view of the result on long time existence of solutions to two-dimensional 
incompressible nonhomogeneous Navier-Stokes equations (see [AKM, Ch. 3]) 
we could expect that smallness of $\varrho_{0,x}$ can be replaced by 
smallness of $\varrho_{0,x_3}$ only.

\noindent
However, up to now, we do not know how to do it.

\noindent
One could expect that looking for solutions close to two-dimensional solutions 
is nothing to do comparing with [AKM, Ch. 3]. But it is totally not true 
because we need three-dimensional imbeddings, solvability of three-dimensional 
problems and apply the three-dimensional technique of Sobolev and 
Sobolev-Slobodetski spaces.

Moreover, we have to mention that many techniques used in this paper were 
developed in [Z2, Z3, Z4, RZ] in the case of a constant density.

The next step in our considerations will be a global existence result which 
can be proved by extending [Z4, NZ1] to the nonhomogeneous fluids.

Finally we expect an existence of global attractor by applying the technique 
of [NZ2].

Many results on existence and estimates of weak solutions to nonhomogeneous 
incompressible Navier-Stokes equations can be found in [P].

\section{2. Notation}

We use isotropic and anisotropic Lebesgue spaces $L_p(Q)$, $Q\in\{\Omega^T,\break
S^T,\Omega,S\}$, $p\in[1,\infty]$; $L_q(0,T;L_p(Q))$, 
$Q\in\{\Omega,S\}$, $p,q\in[1,\infty]$; isotropic and anisotropic 
Sobolev spaces with the norms
$$
\|u\|_{W_p^s(Q)}=\bigg(\sum_{|\alpha|\le s}\intop_Q
|D_x^\alpha u|^pdx\bigg)^{1/p}
$$
and
$$
\|u\|_{W_p^{s,s/2}(Q^T)}=\bigg(\sum_{|\alpha|+2a\le s}\intop_{Q^T}
|D_x^\alpha\partial_t^au|^pdxdt\bigg)^{1/p},\quad s\rm\ \ even,
$$
where $D_x^\alpha=\partial_{x_1}^{\alpha_1}\partial_{x_2}^{\alpha_2}
\partial_{x_3}^{\alpha_3}$, $|\alpha|=\alpha_1+\alpha_2+\alpha_3$, 
$a,\alpha_i\in\N\cup\{0\}\equiv\N_0$, $i=1,2,3$, $Q\in\{\Omega,S\}$, $s\in\N$, 
$p\in[1,\infty]$.

\noindent
In the case $p=2$ we use the notation
$$
H^s(Q)=W_2^s(Q),\quad H^{s,s/2}(Q^T)=W_2^{s,s/2}(Q^T),\quad 
Q\in\{\Omega,S\}.
$$
Moreover, $L_2(Q)=H^0(Q)$, $L_p(Q)=W_p^0(Q)$, $L_p(Q^T)=W_p^{0,0}(Q^T)$.

\noindent
Next we introduce a space natural for examining weak solutions to the 
Navier-Stokes and parabolic equations
$$\eqal{
\|u\|_{V_2^k(\Omega^T)}&=\bigg\{u:\ \|u\|_{V_2^k(\Omega^T)}=
{\rm ess}\sup_{t\in[0,T]}\|u(t)\|_{H^k(\Omega)}\cr
&\quad+\bigg(\intop_0^T\|\nabla u\|_{H^k(\Omega)}^2dt\bigg)^{1/2}<
\infty\bigg\},\quad k\in\N_0.\cr}
$$
In the case of noneven $s$ spaces $W_p^s(Q)$ and $W_p^{s,s/2}(Q^T)$ are 
defined as sets of functions with the finite norms, respectively,
$$\eqal{
\|u\|_{W_p^s(Q)}&=\bigg(\sum_{|\alpha|\le[s]}\intop_Q|D_x^\alpha u(x)|^pdx\cr
&\quad+\sum_{|\alpha|=[s]}\intop_Q\intop_Q
{|D_x^\alpha u(x)-D_{x'}^\alpha u(x')|^p\over|x-x'|^{n+p(s-[s])}}
dxdx'\bigg)^{1/p},\cr}
$$
and
$$\eqal{
&\|u\|_{W_p^{s,s/2}(Q^T)}=\bigg(\sum_{|\alpha|+2a\le[s]}\intop_{Q^T}
|D_x^\alpha\partial_t^au|^pdxdt\cr
&\quad+\sum_{|\alpha|+2a=[s]}\intop_0^T\intop_Q\intop_Q
{|D_x^\alpha\partial_t^au(x,t)-D_{x'}^\alpha\partial_t^au(x',t)|^p\over
|x-x'|^{n+p(s-[s])}}dxdx'dt\cr
&\quad+\sum_{|\alpha|+2a=[s]}\intop_Q\intop_0^T\intop_0^T
{|D_x^\alpha\partial_t^au(x,t)-D_x^\alpha\partial_{t'}^au(x,t')|^p\over
|t-t'|^{1+p(s/2-[s/2])}}dxdtdt'\bigg)^{1/p},\cr}
$$
where $Q\subset\R^n$, $[s]$ the integer parts of $s$.

\noindent
In the case where either $s=[s]$ or ${s\over2}=\left[{s\over2}\right]$ the 
corresponding fractional derivatives vanish.
For $Q=S$ the above norms are defined by applying a partition of unity.

Theorems of imbedding and interpolation for above spaces can be found in [BIN].

By $\mathop C\limits^\krop\null^\alpha(\Omega^T)$ we denote a space of 
functions with the finite seminorm
$$
\|u\|_{\mathop C\limits^\krop\null^\alpha(\Omega^T)}=\sup_{x,x'\in\Omega}
\sup_{t,t'\in(0,T)}{|u(x,t)-u(x',t')|\over|x-x'|^\alpha+|t-t'|^\alpha},
$$
where $\alpha\in(0,1)$.

\noindent
By $c$ we denote a generic constant which changes its value from formula to 
formula. In general $c$ depends on constants of imbeddings, on functions 
describing the boundary, but it does not depend on data. By $\varphi$ we 
denote a generic function which depends on data, changes its form from 
formula to formula and is always positive increasing function of its arguments.

The dependence of $\varphi$ on data will be always expressed explicitly.

To simplify presentation we use the notation
$$
\varrho_x^*=\|\varrho_x\|_{L_\infty(\Omega^T)}.
$$
Let us consider the Stokes system
$$\eqal{
&v_t-\divv\T(v,p)=f'\quad &{\rm in}\ \ \Omega^T,\cr
&\divv v=0\quad &{\rm in}\ \ \Omega^T,\cr
&v\cdot\bar n=0\quad &{\rm on}\ \ S^T,\cr
&\nu\bar n\cdot\D(v)\cdot\bar\tau_\alpha+\delta_{i1}\gamma v\cdot
\bar\tau_\alpha=h_\alpha,\ \ \alpha=1,2,\quad &{\rm on}\ \ S_i^T,\ \ i=1,2,\cr
&v|_{t=0}=v_0\quad &{\rm in}\ \ \Omega.\cr}
\leqno(2.1)
$$
From [Z3] we have

\proclaim Lemma 2.1. 
Assume that $f'\in W_2^{s,s/2}(\Omega^T)$, 
$h_\alpha\in W_2^{s+1/2,s/2+1/4}(S^T)$, $\alpha=1,2$, $s\in\R_+\cup\{0\}$, 
$v_0\in W_2^{s+1}(\Omega)$, $S\in C^{[s]+3}$, where $[s]$ is the integer 
part of $s$.
Then there exists a unique solution to problem (2.1) such that 
$v\in W_2^{s+2,s/2+1}(\Omega^T)$, $\nabla p\in W_2^{s,s/2}(\Omega^T)$ and
$$\eqal{
&\|v\|_{W_2^{s+2,s/2+1}(\Omega^T)}+\|\nabla p\|_{W_2^{s,s/2}(\Omega^T)}\cr
&\le c\Big(\|f'\|_{W_2^{s,s/2}(\Omega^T)}+\|v_0\|_{W_2^{s+1}(\Omega)}+
\sum_{\alpha=1}^2\|h_\alpha\|_{W_2^{s+1/2,s/2+1/4}(S^T)}\Big).\cr}
\leqno(2.2)
$$

\noindent
After small modifications of the proof from [A1] we have

\proclaim Lemma 2.2. 
Assume that $f'\in L_r(\Omega^T)$, $h_\alpha\in W_r^{1-1/r,1/2-1/2r}(S^T)$, 
$r\in(1,\infty)$, $S\in C^2$, $v_0\in W_r^{2-2/r}(\Omega)$.\\
Then there exists a solution to problem (2.1) such that 
$v\in W_r^{2,1}(\Omega^T)$, $\nabla p\in L_r(\Omega^T)$ and
$$\eqal{
&\|v\|_{W_r^{2,1}(\Omega^T)}+\|\nabla p\|_{L_r(\Omega^T)}\le c
\bigg(\|f'\|_{L_r(\Omega^T)}\cr
&\quad+\sum_{\alpha=1}^2\|h_\alpha\|_{W_r^{1-1/r,1/2-1/2r}(S^T)}+
\|v_0\|_{W_r^{2-2/r}(\Omega)}\bigg).\cr}
\leqno(2.3)
$$

\noindent
Let us consider the problem
$$\eqal{
&\varrho v_t-\divv\T(v,p)=\varrho f'\quad &{\rm in}\ \ \Omega^T,\cr
&\divv v=0\quad &{\rm in}\ \ \Omega^T,\cr
&v\cdot\bar n=0\quad &{\rm on}\ \ S^T,\cr
&\nu\bar n\cdot\D(v)\cdot\bar\tau_\alpha+\delta_{i1}\gamma v\cdot
\bar\tau_\alpha=0,\ \ \alpha=1,2,\quad &{\rm on}\ \ S^T,\ i=1,2,\cr
&v|_{t=0}=v_0\quad &{\rm in}\ \ \Omega.\cr}
\leqno(2.4)
$$

\noindent
Using a partition of unity, next the result from [A1] in the half-space and 
finally applying a perturbation argument we have

\proclaim Lemma 2.3. 
Let the assumptions of Lemma 2.2 hold. Let 
$\varrho\in\mathop C\limits^\krop\null^\alpha(\Omega^T)$, $\alpha\in(0,1)$, 
$\varrho,1/\varrho\in L_\infty(\Omega^T)$, 
$\nabla\varrho\in L_\infty(\Omega^T)$ and there exist positive constants 
$0<\varrho_*<\varrho^*$ such that $\varrho_*\le\varrho\le\varrho^*$. Let 
$v\in L_r(\Omega^T)$, $p\in L_r(\Omega^T)$. Then for solutions to problem 
(2.4) the following inequality holds
$$\eqal{
&\|v\|_{W_r^{2,1}(\Omega^T)}+\|\nabla p\|_{L_r(\Omega^T)}\le\varphi
(\varrho_*,\varrho^*,\|\nabla\varrho\|_{L_\infty(\Omega^T)})
[\|v\|_{L_r(\Omega^T)}\cr
&\quad+\|p\|_{L_r(\Omega^T)}+\|f'\|_{L_r(\Omega^T)}+
\|v_0\|_{W_r^{2-2/r}(\Omega)}].\cr}
\leqno(2.5)
$$

\noindent
By a partition of unity and the result from [Z3] in the half-space 
a~perturbation argument implies

\proclaim Lemma 2.4. 
Let the assumptions of Lemma 2.1 be satisfied. Let \\
$\varrho,{1\over\varrho},\varrho_x,\varrho_t\in L_\infty(\Omega^T)$ and let 
there exist positive constants $0<\varrho_*<\varrho^*$ such that 
$\varrho_*\le\varrho\le\varrho^*$. Let $s\in\R_+\cup\{0\}$ and let 
$\varrho\in W_\infty^{s,s/2}(\Omega^T)$. Let $v\in L_2(\Omega^T)$ and 
$p\in W_2^{s,s/2}(\Omega^T)$.
Then there exists a solution to problem (2.4) such that 
$v\in W_2^{s+2,s/2+1}(\Omega^T)$, $\nabla p\in W_2^{s,s/2}(\Omega^T)$, 
$s\in(0,1)$ and
$$\eqal{
&\|v\|_{W_2^{2+s,1+s/2}(\Omega^T)}+\|\nabla p\|_{W_2^{s,s/2}(\Omega^T)}\cr
&\le\varphi(\varrho_*,\varrho^*)(1+\|\partial_x\varrho\|_{L_\infty(\Omega^T)}+
\|\partial_t\varrho\|_{L_\infty(\Omega^T)})\cr
&\quad\cdot[\|v\|_{L_2(\Omega^T)}+\|p\|_{H^{s,s/2}(\Omega^T)}+
\|f'\|_{H^{s,s/2}(\Omega^T)}+\|v_0\|_{H^{1+s}(\Omega)}],\cr}
\leqno(2.6)
$$
where $\varphi$ is the generic function.

\section{3. Auxiliary results}

This section is devoted to obtain some a priori estimates for solutions to 
problem (1.1). Therefore, we assume existence of such solutions to (1.1) that 
the derived estimates can be satisfied.

First we introduce weak solutions

\Def{3.1.} 
By a weak solution to problem (1.1) we mean $v\in V_2^0(\Omega^T)$, 
$\varrho\in L_\infty(\Omega^T)$ such that $\divv v=0$ 
$v\cdot\bar n|_S=0$ and satisfying the integral identities
$$\eqal{
&\intop_{\Omega^T}[-\varrho v\phi_t-\varrho v\otimes v\cdot\nabla\phi+
{\nu\over2}\D(v)\cdot\D(\phi)]dxdt\cr
&\quad+\gamma\intop_{S_1}v\cdot\bar\tau_\alpha\phi\cdot\bar\tau_\alpha dS_1
-\intop_\Omega\varrho_0v_0\phi|_{t=0}dx=\intop_{\Omega^T}\varrho f\phi dx
dt,\cr
&\intop_{\Omega^T}[\varrho\psi_{,t}+\varrho v\cdot\nabla\psi]dxdt+
\intop_\Omega\varrho_0\psi|_{t=0}dx=0,\cr}
\leqno(3.1)
$$
for any $\phi,\psi\in W_{5/2}^{1,1}(\Omega^T)$ such that 
$\phi\cdot\bar n|_S=0$, $\divv\phi=0$, $\phi(T)=0$, $\psi(T)=0$ and the 
summation convenction over the repeated indices is assumed.

\noindent
We need the Korn inequality

\proclaim Lemma 3.2. (see [Z1]) 
Assume that
$$
E_\Omega(v)=\|\D(v)\|_{L_2(\Omega)}^2<\infty,\quad v\cdot\bar n|_S=0,\quad
\divv v=0.
\leqno(3.2)
$$
If $\Omega$ is not axially symmetric there exists a constant $c_1$ such that
$$
\|v\|_{H^1(\Omega)}^2\le c_1E_\Omega(v).
\leqno(3.3)
$$
If $\Omega$ is axially symmetric, $\eta=(-x_2,x_1,0)$, 
$\alpha=\intop_\Omega v\cdot\eta dx$ is bounded then there exists a constant 
$c_2$ such that
$$
\|v\|_{H^1(\Omega)}^2\le c_2\bigg(E_\Omega(v)+
\bigg|\intop_\Omega v\cdot\eta dx\bigg|^2\bigg).
\leqno(3.4)
$$

\noindent
Let us consider the problem
$$\eqal{
&\varrho_t+v\cdot\nabla\varrho=0\quad &{\sl in}\ \ \Omega^T,\cr
&\varrho|_{t=0}=\varrho_0\quad &{\sl in}\ \ \Omega.\cr}
\leqno(3.5)
$$

\proclaim Lemma 3.3. 
Let $\varrho_0\in L_p(\Omega)$, $p\in\R\setminus\{0\}$, $\divv v=0$, 
$v\cdot\bar n|_S=0$. Assume that there exists a sufficiently regular solution 
to problem (3.5). Then the a~priori equality holds
$$
\|\varrho\|_{L_\infty(0,T;L_p(\Omega))}=\|\varrho_0\|_{L_p(\Omega)}.
\leqno(3.6)
$$

\Proof 
Let $p\le1$ and $\varrho\in C^1(\Omega^T)$. Multiplying (3.5) by 
$\varrho|\varrho|^{p-2}$ and integrating over $\Omega^T$ yields, by the 
density argument, the equality
$$
\|\varrho\|_{C(0,T;L_p(\Omega))}=\|\varrho_0\|_{L_p(\Omega)}.
$$
Hence (3.6) holds.

\noindent
For $p<1$, $p\not=0$ we assume additionally that $\varrho,\varrho_0$ are 
different from zero. Hence, performing the same considerations as in the 
case $p\ge1$, we obtain the same equality as above. Finally, (3.6) also holds 
for $p<1$, $p\not=0$. This concludes the proof.

\Remark{3.4.} 
We can pass with $p$ to $+\infty$ and $-\infty$ in (3.6). Let $\varrho_*$, 
$\varrho^*$ be positive constants such that
$$
\varrho_*\le\varrho_0\le\varrho^*.
\leqno(3.7)
$$
Then (3.6) implies
$$
\varrho_*\le\varrho(x,t)\le\varrho^*.
\leqno(3.8)
$$
Next we formulate a result about weak solutions

\proclaim Lemma 3.5. 
Assume that $\Omega$ is not axially symmetric. Assume that 
$f\in L_1(0,T;L_2(\Omega))$, $v_0\in L_2(\Omega)$, 
$\varrho_0,{1\over\varrho_0}\in L_\infty(\Omega)$. Assume that there exist 
constants $\varrho_*$, $\varrho^*$ described in Remark 3.4 and (3.7) holds.
Then a weak solution to problem (1.1) is such that 
$v\in V_2^0(\Omega^T)$ and $\varrho_*\le\varrho(x,t)\le\varrho^*$. Moreover, 
we have the a priori estimates
$$
\|v(t)\|_{L_2(\Omega)}\le\bigg({\varrho^*\over\varrho_*}\bigg)^{1/2}d_1=d_2,
\quad t\le T,
\leqno(3.9)
$$
$$
\|v\|_{V_2^0(\Omega^t)}=\bigg({\varrho^*\over c_3}\bigg)^{1/2}d_1\equiv d_3,
\quad t\le T,
\leqno(3.10)
$$
where $c_3=\min\left({1\over2}\varrho_*,\nu c_1\right)$ and
$$
d_1(T)=\|f\|_{L_1(0,T;L_2(\Omega))}+\|v_0\|_{L_2(\Omega)}.
\leqno(3.11)
$$

\Proof 
Assume that we have the existence of sufficiently regular solutions to (1.1). 
Multiplying $(1.1)_1$ by $v$, integrating over $\Omega$ and using 
$(1.1)_{2,3,4}$ yields
$$
{1\over2}{d\over dt}\|\sqrt{\varrho}v\|_{L_2(\Omega)}^2+\gamma
\|\D(v)\|_{L_2(\Omega)}^2+\gamma\intop_{S_1}|v\cdot\bar\tau_\alpha|^2dS_1=
\intop_\Omega\varrho f\cdot vdx.
\leqno(3.12)
$$
Omitting the second and the third terms on the l.h.s. and applying the 
H\"older inequality to the r.h.s. implies
$$
{1\over2}{d\over dt}\|\sqrt{\varrho}v\|_{L_2(\Omega)}^2\le
\|\sqrt{\varrho}f\|_{L_2(\Omega)}\|\sqrt{\varrho}v\|_{L_2(\Omega)}.
$$
Hence we get
$$
{d\over dt}\|\sqrt{\varrho}v\|_{L_2(\Omega)}\le
\|\sqrt{\varrho}f\|_{L_2(\Omega)}.
\leqno(3.13)
$$
Integrating (3.13) with respect to time yields
$$
\|\sqrt{\varrho}v(t)\|_{L_2(\Omega)}\le
\|\sqrt{\varrho}f\|_{L_1(0,T;L_2(\Omega))}+
\|\sqrt{\varrho_0}v_0\|_{L_2(\Omega)},
\leqno(3.14)
$$
where $t\le T$.

\noindent
Using (3.7) and (3.8) in (3.14) gives (3.9).

\noindent
Integrating (3.12) with respect to time and using (3.3) we obtain
$$\eqal{
&{1\over2}\|\sqrt{\varrho}v\|_{L_2(\Omega)}^2+\nu c_1\intop_0^t
\|v\|_{H^1(\Omega)}^2dt'\le\|\sqrt{\varrho}f\|_{L_1(0,t;L_2(\Omega))}\cr
&\quad\cdot{\rm ess}\sup_{t'\le t}\|\sqrt{\varrho}v\|_{L_2(\Omega)}+{1\over2}
\|\sqrt{\varrho_0}v_0\|_{L_2(\Omega)}^2.\cr}
\leqno(3.15)
$$
In view of (3.14), (3.7) and (3.8) we have (3.10). This concludes the proof.

\noindent
To prove the existence with large data we follow the ideas developed in 
[RZ, Z2,Z4]. therefore we introduce the quantities
$$
h=v_{,x_3},\quad q=p_{,x_3},\quad g=f_{,x_3},\quad \chi=(\rot v)_3,\quad 
F=(\rot f)_3.
\leqno(3.16)
$$

\proclaim Lemma 3.6. 
Let $v$, $\varrho$ be given. Then $(h,q)$ is a solution to the problem
$$\eqal{
&\varrho h_{,t}-\divv\T(h,q)=-\varrho(v\cdot\nabla h+h\cdot\nabla v-g)\quad\cr
&\quad-\varrho_{,x_3}(v_{,t}+v\cdot\nabla v-f)\quad &{\sl in}\ \ \Omega^T,\cr
&\divv h=0\quad &{\rm \sl in}\ \ \Omega^T,\cr
&h\cdot\bar n=0,\ \ 
\bar n\cdot\T(h,q)\cdot\bar\tau_\alpha+\gamma h\cdot\bar\tau_\alpha=0,\ \ 
\alpha=1,2,\quad &{\sl on}\ \ S_1^T,\cr
&h_i=0,\ \ i=1,2,\ \ h_{3,x_3}=0\quad &{\sl on}\ \ S_2^T,\cr
&q|_{S_2}=\varrho f_3\quad &{\sl on}\ \ S_2^T,\cr
&h|_{t=0}=h(0)\quad &{\sl in}\ \ \Omega.\cr}
\leqno(3.17)
$$

\Proof 
$(3.17)_{1,2,3,6}$ follow directly from $(1.1)_{1,2,4,5,6}$ by 
differentiation with respect to $x_3$. Similarly as in [Z2] we show the 
boundary condition $(3.17)_4$. This ends the proof.

\noindent
To formulate problem for $\chi$ we introduce
$$\eqal{
&\bar n|_{S_1}={\nabla\psi\over|\nabla\psi|}={1\over|\nabla\psi|}
(\psi_{,x_1}\psi_{,x_2},0),\cr
&\bar\tau_1|_{S_1}={\nabla^\perp\psi\over|\nabla\psi|}=
{1\over|\nabla\psi|}(-\psi_{,x_2},\psi_{,x_1},0),\cr
&\bar\tau_2|_{S_1}=(0,0,1),\cr
&\bar n|_{S_2}=(0,0,1),\quad \bar\tau_1|_{S_2}=(1,0,0),\quad
\bar\tau_2|_{S_2}=(0,1,0).\cr}
\leqno(3.18)
$$

\proclaim Lemma 3.7. 
Let $\varrho$, $v$, $h$ be given. Then $\chi$ is a solution to the problem
$$\eqal{
&\varrho(\chi_{,t}+v\cdot\nabla\chi)-\nu\Delta\chi=\varrho(F+\chi h_3-
v_{3,x_1}h_2+v_{3,x_2}h_1)\quad\cr
&\quad+\varrho_{,x_1}(v_{2,t}+v\cdot\nabla v_2+f_2)-\varrho_{,x_2}
(v_{1,t}+v\cdot\nabla v_1+f_1)\quad &{\sl in}\ \ \Omega^T,\cr
&\chi=v_i(n_{i,x_j}\tau_{1j}+\tau_{1i,x_j}n_j)+v\cdot\bar\tau_1
(\tau_{12,x_1}-\tau_{11,x_2})\cr
&\quad+{\gamma\over\nu}v_j\tau_{1j}\equiv\chi_*\quad &{\sl on}\ \ S_1^T,\cr
&\chi_{,x_3}=0\quad &{\sl on}\ \ S_2^T,\cr
&\chi|_{t=0}=\chi(0)\quad &{\sl in}\ \ \Omega,\cr}
\leqno(3.19)
$$
where the summation convention over the repeated indices is assumed.

\Proof  $(3.19)_1$ follows from applying two-dimensional rot operator to the 
first two equations of $(1.1)_1$. The boundary condition $(3.19)_2$ was 
proved in [Z2]. This ends the proof.
\goodbreak

To apply the energy type method to problem (3.19) we need

\proclaim Lemma 3.8. 
Let $\chi$ satisfy (3.19). Let $\tilde\chi$ be a solution to the problem
$$\eqal{
&\varrho\tilde\chi_{,t}-\nu\Delta\tilde\chi=0\quad &{\sl in}\ \ \Omega^T,\cr
&\tilde\chi=\chi_*\quad &{\sl on}\ \ S_1^T,\cr
&\tilde\chi_{,x_3}=0\quad &{\sl on}\ \ S_2^T,\cr
&\tilde\chi|_{t=0}=0\quad &{\sl in}\ \ \Omega.\cr}
\leqno(3.20)
$$
Then the function $\chi'=\chi-\tilde\chi$ satisfies
$$\eqal{
&\varrho(\chi'_t+v\cdot\nabla\chi')-\nu\Delta\chi'\cr
&\quad=\varrho(F+\chi h_3-v_{3,x_1}h_2+v_{3,x_2}h_1-v\cdot\nabla\tilde\chi)\cr
&\qquad+\varrho_{,x_1}(f_2+v_{2,t}+v\cdot\nabla v_2)\cr
&\qquad-\varrho_{,x_2}(f_1+v_{1,t}+v\cdot\nabla v_1)
\qquad\qquad\qquad\qquad {\sl in}\ \ \Omega^T,\cr
&\chi'|_{S_1}=0,\cr
&\chi'_{,x_3}|_{S_2}=0,\cr
&\chi'|_{t=0}=\chi(0).\cr}
\leqno(3.21)
$$

\noindent
Let us consider problem (3.5). Then we have

\proclaim Lemma 3.9. 
Assume that $v\in W_2^{s,s/2}(\Omega^T)$, $s>{5\over2}$. Let
$$
X_1=\|\varrho_x(0)\|_{L_\infty(\Omega)}+\|\varrho_t(0)\|_{L_\infty(\Omega)}.
\leqno(3.22)
$$
Then the following a priori estimate is valid
$$\eqal{
&\|\varrho_x(t)\|_{L_p(\Omega)}+\|\varrho_t(t)\|_{L_p(\Omega)}+
\|\varrho\|_{\mathop C\limits^\krop\null^\alpha(\Omega^T)}
\le\varphi(\|v\|_{W_2^{s,s/2}(\Omega^T)})X_1,\cr%
&\alpha\le1-1/p,\quad p\le\infty,\cr}
\leqno(3.23)
$$
where
$$
\|\varrho\|_{\mathop C\limits^\krop\null^\alpha(\Omega^T)}=
\sup_{x,x'\in\Omega\atop t,t'\in(0,T)}
{|\varrho(x,t)-\varrho(x',t')|\over|x-x'|^\alpha+|t-t'|^\alpha},
\leqno(3.24)
$$
where $|x-x'|=\sum_{i=1}^3|x_i-x'_i|$.

\Proof 
For solutions to problem (3.5) we have
$$
{d\over dt}\|\varrho_x\|_{L_p(\Omega)}\le\|v_x\|_{L_\infty(\Omega)}
\|\varrho_x\|_{L_p(\Omega)},
$$
for any $p>1$.

\noindent
Integrating with respect to time yields
$$
\|\varrho_x(t)\|_{L_p(\Omega)}\le\exp\bigg(\intop_0^t
\|v_x(t')\|_{L_\infty(\Omega)}dt'\bigg)\|\varrho_x(0)\|_{L_p(\Omega)},
\leqno(3.25)
$$
where $p\in[1,\infty]$.

\noindent
From $(3.5)_1$ we obtain
$$
\|\varrho_t(t)\|_{L_p(\Omega)}\le\|v\|_{L_\infty(\Omega^T)}
\|\varrho_x\|_{L_p(\Omega)}.
\leqno(3.26)
$$
Let us consider the expression
$$\eqal{
&\sup_{x,x'\in\Omega\atop t,t'\in[0,T]}
{|\varrho(x,t)-\varrho(x',t')|\over
\sum_{i=1}^3|x_i-x'_i|^{1-1/p}+|t-t'|^{1-1/p}}\cr
&\le\sup_{x,x'\in\Omega,\atop t,t'\in[0,T]}\bigg[
{|\varrho(x_1,x_2,x_3,t|-\varrho(x'_1,x_2,x_3,t)|\over|x_1-x'_1|^{1-1/p}}\cr
&\quad+{|\varrho(x'_1,x_2,x_3,t)-\varrho(x'_1,x'_2,x_3,t)|
\over|x_2-x'_2|^{1-1/p}}\cr
&\quad+{|\varrho(x'_1,x'_2,x_3,t)-\varrho(x'_1,x'_2,x'_3,t)|
\over|x_3-x'_3|^{1-1/p}}
+{|\varrho(x',t)-\varrho(x',t')|\over|t-t'|^{1-1/p}}\bigg]\cr
&\le\sup_{x,x'\in\Omega\atop t,t'\in[0,T]}\bigg[\sum_{i=1}^3
\bigg(\intop_{x'_i}^{x_i}|\varrho_{,x_i}|^pdx_i\bigg)^{1/p}+
\bigg(\intop_{t'}^t|\varrho_{,t}|^pdt\bigg)^{1/p}\bigg]\cr
&\le c(\|\varrho_x\|_{L_\infty(\Omega^T)}+\|\varrho_t\|_{L_\infty(\Omega^T)}).
\cr}
\leqno(3.27)
$$
Using the imbedding
$$
\|v\|_{L_\infty(\Omega^T)}+\|v_x\|_{L_2(0,T;L_\infty(\Omega))}\le c
\|v\|_{W_2^{s,s/2}(\Omega^T)}
\leqno(3.28)
$$
for $s>{5\over2}$, we obtain from (3.25)--(3.27) estimate (3.23). This 
concludes the proof.

\proclaim Lemma 3.10. 
Let $\varrho$ and $v$ be given and sufficiently regular. 
Assume also that vectors $\bar n$, $\bar\tau_\alpha$, $\alpha=1,2$, are 
defined in a neighbourhood of $S$ and $a_{\alpha\beta}$, $a$ depend on 
$D_x^\sigma\bar n$, $D_x^\sigma\bar\tau_\alpha$, $\alpha=1,2$, $\sigma\le2$.
Then $p$ is a solution to the problem
$$\eqal{
&\Delta p=-\nabla\varrho\cdot v_t-\divv(\varrho v\cdot\nabla v)+
\divv(\varrho f)\quad {\sl in}\ \ \Omega,\cr
&{\partial p\over\partial n}\bigg|_{S_j}=
(\varrho v_iv_jn_{j,x_i}+\nu a_{j\alpha\beta}v_{\tau_\alpha,\tau_\beta}+
\nu a_j\cdot v+\varrho f\cdot\bar n)|_{S_j},\ \ j=1,2,\cr
&\intop_\Omega pdx=0,\cr}
\leqno(3.29)
$$
where the last equation was added to have uniqueness of solutions to (3.29).

\Proof  Applying div to $(1.1)_1$ we get $(3.29)_1$.
Multiplying $(1.1)_1$ by $\bar n$ and projecting on $S$ we obtain the boundary 
condition
$$\eqal{
{\partial p\over\partial n}\bigg|_S&=
(-\varrho\bar n\cdot v_t-\varrho n_iv\cdot\nabla v_i+\nu\bar n\cdot
\Delta v+\varrho f\cdot\bar n)|_S\cr
&=(\varrho v_iv_jn_{j,x_i}+\nu\bar n\cdot\Delta v+\varrho f\cdot\bar n)|_S
\equiv I,\cr}
\leqno(3.30)
$$
where we used that $v\cdot\bar n|_S=0$ and the summation convention over the 
repeated indices is assumed.

\noindent
Now we calculate $\bar n\cdot\Delta v|_S$. Let us introduce the curvilinear 
coordinates $n$, $\tau_\alpha$, $\alpha=1,2$, connected with the orthonormal 
system of vectors $\bar n$, $\bar\tau_\alpha$, $\alpha=1,2$. Then we calculate
$$\eqal{
\Delta v&=v_{,x_ix_i}=(v_{,n}n_{,x_i}+v_{,\tau_\alpha}\tau_{\alpha,x_i})_{,x_i}
\cr
&=v_{,nn}n_{,x_i}n_{,x_i}+2v_{,n\tau_\alpha}n_{,x_i}\tau_{\alpha,x_i}+
v_{,\tau_\alpha\tau_\beta}\tau_{\alpha,x_i}\tau_{\beta,x_i}\cr
&\quad+v_{,n}n_{,x_ix_i}+v_{,\tau_\alpha}\tau_{\alpha,x_ix_i}\equiv J.\cr}
\leqno(3.31)
$$
By the properties of curvilinear coordinates such that $\bar n\|\nabla n$, 
$\bar\tau_\alpha\|\nabla\tau_\alpha$, $\bar n\cdot\bar\tau_\alpha=0$, 
$\bar\tau_\alpha\cdot\bar\tau_\beta=\delta_{\alpha\beta}$, 
$\alpha,\beta=1,2$, we have
$$
\Delta v=J=v_{,nn}+v_{,\tau_\alpha\tau_\alpha}+v_{,n}\Delta n+v_{,\tau_\alpha}
\Delta\tau_\alpha,
$$
where we used that $n_{,x_i}n_{x_i}=1$, $\tau_{\alpha,x_i}\tau_{\beta,x_i}=
\delta_{\alpha\beta}$, $\alpha,\beta=1,2$.

\noindent
Expressing the equation of continuity in the curvilinear coordinates we have
$$
\divv v\equiv v_{n,n}+v_n\divv\bar n+v_{\tau_\alpha,\tau_\alpha}+
v_{\tau_\alpha}\divv\bar\tau_\alpha=0,
\leqno(3.32)
$$
where $v_n=v\cdot\bar n$, $v_{\tau_\alpha}=v\cdot\bar\tau_\alpha$, 
$\alpha=1,2$.

\noindent
Next we formulate the second boundary condition $(1.1)_5$ in the curvilinear 
coordinates
$$
v_{\tau_\alpha,n}-v_j\tau_{j\alpha,n}-v_in_{i,\tau_\alpha}+\delta_{1j}
\gamma v_{\tau_\alpha}=0\quad {\rm on}\ \ S_j,\ \ j=1,2,\ \alpha=1,2.
\leqno(3.33)
$$
Now, we calculate
$$\eqal{
\bar n\cdot\Delta v|_S&=(v_{n,nn}-\bar n_{,nn}\cdot v-2\bar n_{,n}\cdot v_{,n}+
v_{n,\tau_\alpha\tau_\alpha}\cr
&\quad-\bar n_{,\tau_\alpha\tau_\alpha}\cdot v-2\bar n_{,\tau_\alpha}\cdot
v_{,\tau_\alpha}+\bar n\cdot v_{,n}\Delta n+\bar n\cdot v_{,\tau_\alpha}
\Delta\tau_\alpha)|_S\cr
&=(v_{n,nn}-\bar n_{,nn}\cdot v-2\bar n_{,n}\cdot(v_n\bar n+v_{\tau_\alpha}
\bar\tau_\alpha)_{,n}-\bar n_{,\tau_\alpha\tau_\alpha}\cdot v\cr
&\quad-2\bar n_{,\tau_\alpha}\cdot(v_n\bar n+v_{\tau_\beta}
\bar\tau_\beta)_{,\tau_\alpha}+(v_{n,n}-\bar n_{,n}\cdot v)\Delta n\cr
&\quad-\bar n_{,\tau_\alpha}\cdot v\Delta\tau_\alpha)|_S\cr
&=(v_{n,nn}-\bar n_{,nn}\cdot v-2\bar n_{,n}\cdot\bar\tau_\alpha 
v_{\tau_\alpha,n}-2\bar n_{,n}\cdot\bar\tau_{\alpha,n}v_{\tau_\alpha}\cr
&\quad-\bar n_{,\tau_\alpha\tau_\alpha}\cdot v-2\bar n_{,\tau_\alpha}\cdot
\bar\tau_\beta v_{\tau_\beta,\tau_\alpha}-2\bar n_{,\tau_\alpha}
\bar\tau_{\beta,\tau_\alpha}v_{\tau_\beta}\cr
&\quad+(v_{n,n}-\bar n_{,n}\cdot v)\Delta n-\bar n_{,\tau_\alpha}\cdot v
\Delta\tau_\alpha)|_S.\cr}
\leqno(3.34)
$$
From (3.32) we calculate
$$\eqal{
v_{n,nn}&=-v_{n,n}\divv\bar n-v_{\tau_\alpha,\tau_\alpha n}-
v_{\tau_\alpha,n}\divv\bar\tau_\alpha\cr
&\quad-v_n(\divv\bar n)_{,n}-v_{\tau_\alpha}(\divv\bar\tau_\alpha)_{,n}.\cr}
\leqno(3.35)
$$
Projecting (3.35) on $S$ and using (3.32) and (3.33) we obtain
$$\eqal{
&v_{n,nn}|_{S_j}=(v_{\tau_\alpha,\tau_\alpha}+v_{\tau_\alpha}
\divv \bar\tau_\alpha)\divv\bar n-(v\cdot\bar\tau_{\alpha,n}+v\cdot
\bar n_{,\tau_\alpha}-\delta_{1j}\gamma v_{\tau_\alpha})_{,\tau_\alpha}\cr
&\quad-(v\cdot\bar\tau_{\alpha,n}+v\cdot\bar n_{,\tau_\alpha}-
\delta_{1j}\gamma v_{\tau_\alpha})\divv\bar\tau_\alpha-
v_{\tau_\alpha}(\divv\bar\tau_\alpha)_{,n},\ \ j=1,2.\cr}
\leqno(3.36)
$$
Calculating $v_{n,nn}|_S$ from (3.36), $v_{n,n}|_S$ from (3.32) and 
$v_{\tau_\alpha,n}|_S$ from (3.33) and inserting them into (3.34) we obtain
$$
\bar n\cdot\Delta v|_{S_j}=a_{j\alpha\beta}v_{\tau_\alpha,\tau_\beta}+
a_j\cdot v,\ \ j=1,2,
\leqno(3.37)
$$
where for $S_j\in C^\alpha$ we have that $a_{j\alpha\beta}\in C^{\alpha-2}$, 
$a_j\in C^{\alpha-3}$, $j=1,2$.

\noindent
From (3.30), (3.31) and (3.37) we obtain (3.29). This concludes the proof.

Now we estimate the norms $\|p\|_{L_\sigma(\Omega^t)}$, $\sigma={5\over3},2$. 
For this purpose we examine problem (3.29). Let $G$ be the Green function to 
the Neumann problem (3.29). Then any solution to (3.29) can be expressed in 
the form
$$\eqal{
&p(x,t)=\intop_\Omega G(x,y)[-\nabla\varrho\cdot v_t-\divv
(\varrho v\cdot\nabla v)+\divv(\varrho f)]dx\cr
&\quad-\sum_{j=1}^2\intop_{S_j}G(x,y)[-\varrho v\cdot\nabla v\cdot\bar n+
\nu a_{j\alpha\beta}v_{\tau_\alpha,\tau_\beta}+\nu a_j\cdot v+\varrho f\cdot
\bar n]dS_{jy}.\cr}
\leqno(3.38)
$$
Integrating by parts in the second and the third expressions of the first 
integral and in the second expression of the second integral of (3.38) we get
$$\eqal{
p(x,t)&=\intop_\Omega[G(x,y)(-\nabla\varrho\cdot v_t)+(\varrho v\cdot
\nabla v-\varrho f)\nabla_yG(x,y)]dy\cr
&\quad-\sum_{j=1}^2\intop_{S_j}\{G(x,y)[-\nu a_{j\alpha\beta_{,\tau_\beta}}
v_{\tau_\alpha}+\nu a_j\cdot v]\cr
&\quad-\nu G(x,y)_{,\tau_\beta}a_{j\alpha\beta}v_{\tau_\alpha}\}dS_{jy}.\cr}
\leqno(3.39)
$$

\proclaim Lemma 3.11. 
Assume that $v\in W_{5/3}^{2,1}(\Omega^T)$, 
$\varrho\in L_\infty(0,T;W_\infty^1(\Omega))$, $v\in V_2^0(\Omega^T)$, 
$f\in L_{5/3}(0,T;L_{15/14}(\Omega))$. Then the following inequality holds
$$\eqal{
&\|p\|_{L_{5/3}(\Omega^T)}\le\varepsilon\|v\|_{W_{5/3}^{2,1}(\Omega^T)}+
c\varrho_x^*\|v_t\|_{L_{5/3}(\Omega^T)}\cr
&\quad+c(1/\varepsilon,\varrho^*)d_3+c\varrho^*
\|f\|_{L_{5/3}(0,T;L_{15/14}(\Omega))},\cr}
\leqno(3.40)
$$
where $\varepsilon\in(0,1)$.

\Proof 
By the properties of the Green function we obtain
$$\eqal{
\|p\|_{L_{5/3}(\Omega)}&\le c[\varrho_x^*\|v_t\|_{L_{5/3}(\Omega)}+\varrho^*
\|v\cdot\nabla v\|_{L_{15\over14}(\Omega)}+\varrho^*
\|f\|_{L_{15\over14}(\Omega)}\cr
&\quad+\|v\|_{W_{15/14}^{1-14/15}(S)}]\cr
&\le c[\varrho_x^*\|v_t\|_{L_{5/3}(\Omega)}+\varrho^*
\|v\cdot\nabla v\|_{L_{15/14}(\Omega)}\cr
&\quad+\|v\|_{W_{15/14}^1(\Omega)}+\varrho^*\|f\|_{L_{15/14}(\Omega)}]\cr
&\le c[\varrho_x^*\|v_t\|_{L_{5/3}(\Omega)}+\varrho^*
\|v\|_{L_{30/13}(\Omega)}\|\nabla v\|_{L_2(\Omega)}\cr
&\quad+\|v\|_{W_{15/14}^1(\Omega)}+
\varrho^*\|f\|_{L_{15/14}(\Omega)}].\cr}
$$
Integrating with respect to time we get
$$\eqal{
\|p\|_{L_{5/3}(\Omega^T)}&\le c[\varrho_x^*\|v_t\|_{L_{5/3}(\Omega^T)}+
\varrho^*\|v\|_{L_{10}(0,T;L_{30\over13}(\Omega))}
\|\nabla v\|_{L_2(\Omega^T)}\cr
&\quad+\|v\|_{L_{5/3}(0,T;W_{15/14}^1(\Omega))}+\varrho^*
\|f\|_{L_{5/3}(0,T';L_{15/14}(\Omega))}].\cr}
$$
By certain interpolation and the energy type estimate we obtain
$$\eqal{
\|p\|_{L_{5/3}(\Omega^T)}&\le\varepsilon\|v\|_{W_{5/3}^{2,1}(\Omega^T)}+
c\varrho_x^*\|v_t\|_{L_{5/3}(\Omega^T)}\cr
&\quad+c(1/\varepsilon,\varrho^*)\|v\|_{L_2(0,T;H^1(\Omega))}+c\varrho^*
\|f\|_{L_{5/3}(0T;L_{15/14}(\Omega))}.\cr}
\leqno(3.41)
$$
In view of (3.10) we get (3.40). This concludes the proof.

\noindent
Next, we have

\proclaim Lemma 3.12. 
Assume that $v\in L_\infty(0,T;L_3(\Omega))\cap V_2^0(\Omega^T)$, 
$v_t\in L_2(\Omega^T)$, $f\in L_2(0,T;L_{6/5}(\Omega))$.\\
Then the following inequality is valid
$$\eqal{
&\|p\|_{L_2(\Omega^T)}\le c\varrho_x^*\|v_t\|_{L2(\Omega^T)}+c(\varrho^*
\|v\|_{L_\infty(0,T;L_3(\Omega))}+1)d_3\cr
&\quad+c\varrho^*\|f\|_{L_2(0,T;L_{6/5}(\Omega))}).\cr}
\leqno(3.42)
$$

\Proof 
By the properties of the Green function we have also
$$\eqal{
\|p\|_{L_2(\Omega)}&\le c\varrho_x^*\|v_t\|_{L_2(\Omega)}+c\varrho^*
(\|v\cdot\nabla v\|_{L_{6/5}(\Omega)}+\|f\|_{L_{6/5}(\Omega)})\cr
&\quad+c\|v\|_{W_{6/5}^1(\Omega)}.\cr}
$$

Integrating with respect to time is the $L_2$-norm and using the H\"older 
inequality we obtain
$$\eqal{
\|p\|_{L_2(\Omega^T)}&\le c\varrho_x^*\|v_t\|_{L_2(\Omega^T)}+
c\varrho^*\|v\|_{L_\infty(0,T;L_3(\Omega))}+1)\|v\|_{L_2(0,T;H^1(\Omega))}\cr
&\quad+c\varrho^*\|f\|_{L_2(0,T;L_{6/5}(\Omega))}.\cr}
\leqno(3.41)
$$
Using (3.10) we get (3.42). This concludes the proof.

\section{4. Estimates}

First we obtain an estimate for solutions to problem (3.19).

\proclaim Lemma 4.1. 
Assume that $\varrho\in L_\infty(0,T;W_\infty^1(\Omega))$, 
$\varrho_x^*=\|\varrho_x\|_{L_\infty(\Omega^T)}$, 
$\varrho\in C^\alpha(\Omega^T)$, $v'=(v_1,v_2)$, 
$v'\in L_5(\Omega^T)\cap L_\infty(0,T;L_2(\Omega))\cap W_2^{1,1/2}(\Omega^T)$, 
$v'_t\in L_2(0,T;L_{6/5}(\Omega))$, $\nabla v'\in L_2(0,T;L_3(\Omega))$, 
$h\in L_\infty(0,T;L_3(\Omega))$, $F\in L_2(0,T;L_{6/5}(\Omega))$, 
$f'\in L_2(0,T;L_{6/5}(\Omega))$, $f'=(f_1,f_2)$, $\chi(0)\in L_2(\Omega)$. 
Assume that $v$ is a weak solution to problem (1.1).\\
Then solutions to problem (3.19) satisfy the inequality
$$\eqal{
&\sigma_1\|\chi\|_{V_2^0(\Omega^t)}\le c\varrho^*d_3
(\|h\|_{L_\infty(0,t;L_3(\Omega))}+\|F\|_{L_2(0,t;L_{6/5}(\Omega))})\cr
&\quad+c\varrho_x^*(\|f'\|_{L_2(0,t;L_{6/5}(\Omega))}+
\|v'_t\|_{L_2(0,t;L_{6/5}(\Omega))}\cr
&\quad+d_2\|\nabla v'\|_{L_2(0,t;L_3(\Omega))})
+c\varrho^*\|\chi(0)\|_{L_2(\Omega)}\cr
&\quad+\varphi(\varrho_*,\varrho^*,\varrho_x^*,
\|\varrho\|_{\mathop C\limits^\krop\null^\alpha(\Omega^t)})
(\|v'\|_{L_5(\Omega^t)}+d_2\cr
&\quad+\|v'\|_{W_2^{1,1/2}(\Omega^t)}),\quad t\le T,\cr}
\leqno(4.1)
$$
$\sigma_1=\min\{\varrho_*,\nu\}$, $\varphi$ is an increasing continuous 
positive function.
\goodbreak

\Proof 
Multiplying $(3.21)_1$ by $\chi'$, integrating over $\Omega$, using the 
continuity equation $(1.1)_{2,3}$ and the boundary conditions yields
$$\eqal{
&{d\over dt}\intop_\Omega\varrho\chi^{'2}dx+\nu\intop_\Omega
|\nabla\chi'|^2dx=\intop_\Omega\varrho h_3\chi\chi'dx\cr
&\quad-\intop_\Omega\varrho v\cdot\nabla\tilde\chi\chi'dx+\intop_\Omega
\varrho(F-v_{3,x_1}h_2+v_{3,x_2}h_1)\chi'dx\cr
&\quad+\intop_\Omega[\varrho_{x_1}(f_2+v_{2,t}+v\cdot\nabla v_2)-
\varrho_{,x_2}(f_1+v_{1,t}+v\cdot\nabla v_1)]\chi'dx.\cr}
\leqno(4.2)
$$
Now we estimate the particular terms from the r.h.s. of (4.2). The first 
term we estimate by
$$
\varepsilon\|\chi'\|_{L_6(\Omega)}^2+c(1/\varepsilon)\varrho^{*2}
\|h_3\|_{L_3(\Omega)}^2\|\chi\|_{L_2(\Omega)}^2,
$$
the second we express in the form 
$\intop_\Omega\varrho v\cdot\nabla\chi'\tilde\chi dx+\intop_\Omega v\cdot
\nabla\varrho\tilde\chi\chi'dx$ and estimate by
$$
\varepsilon\|\nabla\chi'\|_{L_2(\Omega)}^2+c(1/\varepsilon)
(\varrho^{*2}+\|\nabla\varrho\|_{L_3(\Omega)}^2
\|v\tilde\chi\|_{L_2(\Omega)}^2),
$$
the third by
$$
\varepsilon\|\chi'\|_{L_6(\Omega)}^2+c(1/\varepsilon)\varrho^{*2}
(\|F\|_{L_{6/5}(\Omega)}^2+\|v_{3,x'}\|_{L_2(\Omega)}^2\|h'\|_{L_3(\Omega)}^2),
$$
where $h'=(h_1,h_2)$ and finally the last by
$$
\varepsilon\|\chi'\|_{L_6(\Omega)}^2+c(1/\varepsilon)\varrho_x^{*2}
(\|f'\|_{L_{6/5}(\Omega)}^2+\|v'_t\|_{L_{6/5}(\Omega)}^2+
\|v\|_{L_2(\Omega)}^2\|\nabla v'\|_{L_3(\Omega)}^2).
$$
Using the above estimates in (4.2), assuming that $\varepsilon$ is 
sufficiently small, integrating the result with respect to time and using 
Lemma 3.5 we obtain
$$\eqal{
&\sigma_1\|\chi'\|_{V_2^0(\Omega^t)}\le c\varrho^{*2}d_3^2
(\|h\|_{L_\infty(0,t;L_3(\Omega))}^2+\|\tilde\chi\|_{L_5(\Omega^t)}^2\cr
&\quad+\|F\|_{L_2(0,t;L_{6/5}(\Omega))}^2)+c\varrho_x^{*2}
(\|f'\|_{L_2(0,t;L_{6/5}(\Omega))}^2+\|v'_t\|_{L_2(0,t;L_{6/5}(\Omega))}^2\cr
&\quad+d_3^2\|\tilde\chi\|_{L_5(\Omega^t)}^2+d_2^2\|\nabla v'\|_{L_2(0,t;L_3(\Omega))}^2)+\varrho^{*2}
\|\chi(0)\|_{L_2(\Omega)}^2.\cr}
$$
In view of the relation between $\chi$ and $\chi'$ we have
$$\eqal{
&\sigma_1\|\chi\|_{V_2^0(\Omega^t)}^2\le c\varrho^{*2}d_3^2
(\|h\|_{L_\infty(0,t;L_3(\Omega))}^2+\|F\|_{L_2(0,t;L_{6/5}(\Omega))}^2)\cr
&\quad+c\varrho_x^{*2}(\|f'\|_{L_2(0,t;L_{6/5}(\Omega))}^2+
\|v'_t\|_{L_2(0,t;L_{6/5}(\Omega))}^2\cr
&\quad+d_2^2\|\nabla v'\|_{L_2(0,t;L_3(\Omega))}^2)+
\varrho^{*2}\|\chi(0)\|_{L_2(\Omega)}^2\cr
&\quad+c(\varrho^{*2}+\varrho_x^{*2})
d_3^2\|\tilde\chi\|_{L_5(\Omega^t)}^2+
\sigma_1\|\tilde\chi\|_{V_2^0(\Omega^t)}^2.\cr}
\leqno(4.3)
$$
For solutions to problem (3.20) we obtain (see [Z6])
$$\eqal{
&\|\tilde\chi\|_{L_5(\Omega^T)}\le\varphi(\varrho_*,\varrho^*,
\|\varrho\|_{\mathop C\limits^\krop\null^\alpha(\Omega^T)})
\|v'\|_{L_5(\Omega^T)},\cr
&\|\tilde\chi\|_{L_\infty(0,T;L_2(\Omega))}\le\varphi(\varrho_*,\varrho^*,
\|\varrho\|_{\mathop C\limits^\krop\null^\alpha(\Omega^T)})
\|v'\|_{L_\infty(0,T;L_2(\Omega))},\cr
&\|\nabla\tilde\chi\|_{L_2(\Omega^T)}\le\varphi(\varrho_*,\varrho^*,
\|\varrho\|_{\mathop C\limits^\krop\null^\alpha(\Omega^T)})
\|v'\|_{W_2^{1,1/2}(\Omega^T)},\cr}
\leqno(4.4)
$$
where $\varphi$ is an increasing continuous positive function.

\noindent
Using (4.4) in (4.3) implies (4.1). This concludes the proof.

\noindent
Next we shall obtain an estimate for solutions to problem (3.17).

\proclaim Lemma 4.2. 
Assume that $v$ is a weak solution to problem (1.1). Assume that 
$g\in L_2(0,T;L_{6/5}(\Omega))$, $h(0)\in L_2(\Omega)$, 
$f_3\in L_2(0,T;L_{4/3}(S_2))$, $f\in L_2(0,T;L_{6/5}(\Omega))$, 
$v_t\in L_2(0,T;L_{6/5}(\Omega))$, $v\in L_\infty(0,T;L_3(\Omega))$, 
$\varrho\in L_\infty(0,T;W_\infty^1(\Omega))$.\\
Assuming additionally that $h\in L_\infty(0,T;L_3(\Omega))$ we obtain
$$\eqal{
&\sigma_1\|h\|_{V_2^0(\Omega^t)}^2\le c\varrho^{*2}
(\|h\|_{L_\infty(0,t;L_3(\Omega))}^2d_3^2\cr
&\quad+\|g\|_{L_2(0,t;L_{6/5}(\Omega))}^2)+c\varrho_x^{*2}
(\|v_t\|_{L_2(0,t;L_{6/5}(\Omega))}^2\cr
&\quad+d_3^2\|v\|_{L_\infty(0,t;L_3(\Omega))}^2+
\|f\|_{L_2(0,t;L_{6/5}(\Omega))}^2)\cr
&\quad+c(\|f_3\|_{L_2(0,t;L_{4/3}(S_2))}^2+\|h(0)\|_{L_2(\Omega)}^2),\quad
t\le T,\cr}
\leqno(4.5)
$$
where $\sigma_1$ and $\varrho_x^*$ are the same as in Lemma 4.1.\\
Replacing the condition $h\in L_\infty(0,T;L_3(\Omega))$ by 
$v\in L_2(0,T;W_3^1(\Omega))$ we have
$$\eqal{
&\sigma_1\|h\|_{V_2^0(\Omega^t)}^2\le c\exp
(\|\nabla v\|_{L_2(0,t;L_3(\Omega))}^2)\cr
&\quad\cdot[\|g\|_{L_2(0,t;L_{6/5}(\Omega))}^2+\|h(0)\|_{L_2(\Omega)}^2+
\varrho_x^{*2}(\|v_t\|_{L_2(0,t;L_{6/5}(\Omega))}^2\cr
&\quad+d_3^2\|v\|_{L_\infty(0,t;L_3(\Omega))}^2+
\|f\|_{L_2(0,t;L_{6/5}(\Omega)}^2)+\|f_3\|_{L_2(0,t;L_{4/3}(S_2))}^2],\cr}
\leqno(4.6)
$$
where $t\le T$.

\Proof 
Multiplying $(3.17)_1$ by $h$, integrating the result over $\Omega$ and 
using the continuity equation $(1.1)_3$ we obtain
$$\eqal{
&{1\over2}{d\over dt}\intop_\Omega\varrho h^2dx+{\nu\over2}
\|\D(h)\|_{L_2(\Omega)}^2-\intop_S\bar n\cdot\T(h,q)\cdot hdS\cr
&=-\intop_\Omega\varrho h\cdot\nabla v\cdot hdx+\intop_\Omega\varrho g\cdot
hdx-\intop_\Omega\varrho_{,x_3}(v_t+v\cdot\nabla v-f)\cdot hdx.\cr}
\leqno(4.7)
$$
The boundary term in (4.7) equals
$$
{\gamma\over2}\intop_{S_1}|h\cdot\bar\tau_\alpha|^2dS_1-
\intop_{S_2}qh_3dS_2\equiv I
$$
so the second term in $I$ is estimated by
$$
\varepsilon\|h\|_{H^1(\Omega)}^2+c(1/\varepsilon)
\|f_3\|_{L_{4\over3}(S_2)}^2.
$$
The first term on the r.h.s. of (4.7) we estimate in two different ways. 
Either by
$$
\varepsilon\|h\|_{L_6(\Omega)}^2+c(1/\varepsilon)\varrho^{*2}
\|\nabla v\|_{L_2(\Omega)}^2\|h\|_{L_3(\Omega)}^2
\leqno(4.8)
$$
or by
$$
\varepsilon\|h\|_{L_6(\Omega)}^2+c(1/\varepsilon)\varrho^*\intop_\Omega
\varrho h^2dx\|\nabla v\|_{L_3(\Omega)}^2.
\leqno(4.9)
$$
The second term on the r.h.s. of (4.7) we estimate by
$$
\varepsilon\|h\|_{L_6(\Omega)}^2+c(1/\varepsilon)\varrho^{*2}
\|g\|_{L_{6/5}(\Omega)}^2
$$
and the last by
$$
\varepsilon\|h\|_{L_6(\Omega)}^2+c(1/\varepsilon)\varrho_x^{*2}
(\|v_t\|_{L_{6/5}(\Omega)}^2+\|v\cdot\nabla v\|_{L_{6/5}(\Omega)}^2+
\|f\|_{L_{6/5}(\Omega)}^2).
$$
Using he above estimates in (4.7), assuming that $\varepsilon$ is sufficiently 
small, using Lemma 3.5 and integrating with respect to time we obtain (4.5) 
in the case (4.8) and (4.6) for (4.9). Let us mention that the time integral 
of the first term in $I$ is deleted. This concludes the proof.

\noindent
Let us consider the elliptic problem
$$\eqal{
&v_{1,x_2}-v_{2,x_1}=\chi\quad &{\rm in}\ \ \Omega',\cr
&v_{1,x_1}+v_{2,x_2}=-h_3\quad &{\rm in}\ \ \Omega',\cr
&v'\cdot\bar n|_{S'_1}=0.\cr}
\leqno(4.10)
$$
Let $P$ be a plane perpendicular to the axis of the cylinder. Then 
$\Omega'=\Omega\cap P$, $S'_1=S_1\cap P$.

\noindent
In view of (4.1) and (4.5) we obtain for solutions to problem (4.10) the 
inequality
$$\eqal{
&\|v'\|_{V_2^1(\Omega^t)}\le\varphi(\varrho_*,\varrho^*,
\|\varrho\|_{\mathop C\limits^\krop\null^\alpha(\Omega^T)},d_3)
[\|h\|_{L_\infty(0,t;L_3(\Omega))}\cr
&\quad+\|v'\|_{L_5(\Omega^t)}+\|v'\|_{W_2^{1,1/2}(\Omega^t)}+
G_1(t)]+c\varrho_x^*[d_2\|\nabla v'\|_{L_2(0,t;L_3(\Omega))}\cr
&\quad+\|v_t\|_{L_2(0,t;L_{6/5}(\Omega))}+
\|\nabla v'\|_{L_2(0,t;L_3(\Omega))}\cr
&\quad+\|v\|_{L_\infty(0,t;L_3(\Omega))}+\|f\|_{L_2(0,t;L_{6/5}(\Omega))}],
\ \ t\le T,\cr}
\leqno(4.11)
$$
where
$$\eqal{
G_1(t)&=\|F\|_{L_2(0,t;L_{6/5}(\Omega))}+\|g\|_{L_2(0,t;L_{6/5}(\Omega))}\cr
&\quad+\|f_3\|_{L_2(0,t;L_{4/3}(S_2))}+\|\chi(0)\|_{L_2(\Omega)}+
\|h(0)\|_{L_2(\Omega)}.\cr}
\leqno(4.12)
$$
Applying interpolation inequalities in (4.11) (see [BIN, Ch. 3, Sect. 10]) 
implies the inequality
$$\eqal{
\|v'\|_{V_2^1(\Omega^t)}&\le\varphi(\varrho_*,\varrho^*,
\|\varrho\|_{\mathop C\limits^\krop\null^\alpha(\Omega^T)},d_3)
[\|h\|_{L_\infty(0,t;L_3(\Omega))}\cr
&\quad+\|v'\|_{L_2(\Omega;H^{1/2}(0,t))}+d_3+G_1(t)]\cr
&\quad+c\varrho_x^*[\|v_t\|_{L_2(0,t;L_{6/5}(\Omega))}+
\|v\|_{L_\infty(0,t;L_3(\Omega))}\cr
&\quad+\|\nabla v'\|_{L_2(0,t;L_3(\Omega))}+\varphi(\varrho_x^*)d_3\cr
&\quad+\|f\|_{L_2(0,t;L_{6/5}(\Omega))}].\cr}
\leqno(4.13)
$$
In view of (3.23) we have
$$\eqal{
\|v'\|_{V_2^1(\Omega^t)}&\le\varphi(\varrho_*,\varrho^*,\varphi
(T^{1/2}\|v\|_{W_2^{\sigma,\sigma/2}(\Omega^T)})X_1,d_3)\cr
&\quad\cdot[\|h\|_{L_\infty(0,t;L_3(\Omega))}+
\|v'\|_{L_2(\Omega;W_2^{1/2}(0,t))}\cr
&\quad+X_1(\|v_t\|_{L_2(0,t;L_{6/5}(\Omega))}+
\|v\|_{L_\infty(0,t;L_3(\Omega))}\cr
&\quad+\|\nabla v'\|_{L_2(0,t;L_3(\Omega))})+G_2(t)],\quad t\le T,\cr}
\leqno(4.14)
$$
where $\sigma>{5\over2}$,
$$
G_2(t)=G_1(t)+\|f\|_{L_2(0,t;L_{6/5}(\Omega))}+d_3
\leqno(4.15)
$$
and $G_1$ is defined by (4.12) and $X_1$ by (3.22).

\noindent
Now, we consider problem (1.1) in the form
$$\eqal{
&\varrho v_t-\divv\T(v,p)=-\varrho v'\cdot\nabla v-\varrho v_3h+\varrho f,\cr
&\divv v=0\cr
&v\cdot\bar n|_S=0\cr
&\bar n\cdot\T(v,p)\cdot\bar\tau_\alpha+\delta_{1j}\gamma v\cdot
\bar\tau_\alpha|_{S_j}=0,\ \ \alpha=1,2,\ \ j=1,2,\cr
&v|_{t=0}=v(0).\cr}
\leqno(4.16)
$$

\proclaim Lemma 4.3. 
Assume that $f\in L_2(\Omega^T)$, $v_0\in H^1(\Omega)$, 
$v\in W_2^{\sigma,\sigma/2}(\Omega^T)$, $\sigma>{5\over2}$, 
$\varrho_*\le\varrho\le\varrho^*$, $\varrho_x(0)\in L_\infty(\Omega)$, 
$h\in L_\infty(0,T;L_3(\Omega))\cap L_{10\over3}(\Omega^T)$.
Then for solutions to (4.16) the following inequality holds
$$\eqal{
&\|v\|_{W_2^{2,1}(\Omega^T)}+\|\nabla p\|_{L_2(\Omega^T)}\le\varphi
(\varrho_*,\varrho^*,\varphi(T^{1/2}\|v\|_{W_2^{\sigma,\sigma/2}(\Omega^T)})
X_1,d_3)\cr
&\quad\cdot[H^2+X_1^2(\|v_t\|_{L_2(\Omega^T)}^2+
\|v\|_{L_\infty(0,T;L_3(\Omega))}^2+\|\nabla v'\|_{L_2(0,T;L_3(\Omega))}^2)+
G^2]\cr
&\equiv\varphi(H^2+X_1^2V^2+G^2),\cr}
\leqno(4.17)
$$
where $X_1$ is introduced in Lemma 3.9 (see (3.22)),
$$\eqal{
H&=\|h\|_{L_\infty(0,T;L_3(\Omega))}+\|h\|_{L_{10\over3}(\Omega^T)},\cr
G&=\|f\|_{L_2(\Omega^T)}+\|v_0\|_{H^1(\Omega)}+d_3+
\|F\|_{L_2(0,T;L_{6/5}(\Omega))}\cr
&\quad+\|g\|_{L_2(0,T;L_{6/5}(\Omega))}+\|f_3\|_{L_2(0,T;L_{4/3}(S_2))},\cr
V&=\|v_t\|_{L_2(\Omega^T)}+\|v\|_{L_\infty(0,T;L_3(\Omega))}+
\|\nabla v'\|_{L_2(0,T;L_3(\Omega))}\cr}
\leqno(4.18)
$$
and $\varphi$ is a generic function described by the r.h.s. of the above
inequality.

\Proof 
From (2.5) (see Lemma 2.3), energy estimate (3.10), \\
$f'=-\varrho v'\cdot\nabla'v-\varrho v_3h+\varrho f$ we obtain
$$\eqal{
&\|v\|_{W_{5/3}^{2,1}(\Omega^T)}+\|\nabla p\|_{L_{5/3}(\Omega^T)}\cr
&\le\Phi[\|v\|_{L_{5/3}(\Omega^T)}+\|p\|_{L_{5/3}(\Omega^T)}+\varrho^*
\|v'\|_{L_{10}(\Omega^T)}d_3\cr
&\quad+\varrho^*d_3\|h\|_{L_{10\over3}(\Omega^T)}+\varrho^*
\|f\|_{L_{5/3}(\Omega^T)}+\|v_0\|_{W_{5/3}^{4/5}(\Omega)}],\cr}
\leqno(4.19)
$$
where we used that
$$
\|v'\nabla v\|_{L_{5/3}(\Omega^T)}\le\|v'\|_{L_{10}(\Omega^T)}
\|\nabla v\|_{L_2(\Omega^T)},
$$
$$
\|v_3h\|_{L_{5/3}(\Omega^T)}\le\|v_3\|_{L_{10/3}(\Omega^T)}
\|h\|_{L_{10/3}(\Omega^T)}
$$
and we introduced the quantity
$$
\Phi=\varphi(\varrho_*,\varrho^*,\varphi
(T^{1/2}\|v\|_{W_2^{\sigma,\sigma/2}(\Omega^T)})X_1,d_3).
\leqno(4.20)
$$
In view of the imbedding (see [Z3, Lemma 3.7])
$$
\|v'\|_{L_{10}(\Omega^T)}\le c\|v'\|_{V_2^1(\Omega^T)},
\leqno(4.21)
$$
notation (4.18) we obtain, from (4.19) after some interpolations, the 
inequality
$$\eqal{
&\|v\|_{W_{5/3}^{2,1}(\Omega^T)}+\|\nabla p\|_{L_{5/3}(\Omega^T)}\cr
&\le\Phi[\|p\|_{L_{5/3}(\Omega^T)}+H+X_1(\|v_t\|_{L_2(0,T;L_{6/5}(\Omega))}+
\|v\|_{L_\infty(0,T;L_3(\Omega))}\cr
&\quad+\|\nabla v'\|_{L_2(0,T;L_3(S_2))})+G_3],\cr}
\leqno(4.22)
$$
where
$$
G_3=G_2+\|f\|_{L_{5/3}(\Omega^T)}+\|v_0\|_{W_{5/3}^{4/5}(\Omega)}
\leqno(4.23)
$$
and $G_2$ is defined by (4.15).
\goodbreak

\noindent
In view of (3.40) we obtain from (4.22) the inequality
$$\eqal{
&\|v\|_{W_{5/3}^{2,1}(\Omega^T)}+\|\nabla p\|_{L_{5/3}(\Omega^T)}\le
\Phi[H+X_1(\|v_t\|_{L_{5/3}(\Omega^T)}\cr
&\quad+\|v\|_{L_\infty(0,T;L_3(\Omega))}+\|\nabla v'\|_{L_2(0,T;L_3(\Omega))})+
G_3].\cr}
\leqno(4.24)
$$
Using (4.24) in (4.14) yields
$$\eqal{
&\|v'\|_{V_2^1(\Omega^T)}\le\Phi[H+X_1(\|v_t\|_{L_{5/3}(\Omega^T)}\cr
&\quad+\|v\|_{L_\infty(0,T;L_3(\Omega))}+\|\nabla v'\|_{L_2(0,T;L_3(\Omega))})+
G_3].\cr}
\leqno(4.25)
$$
From (2.5) (see Lemma 2.3) we have
$$\eqal{
&\|v\|_{W_2^{2,1}(\Omega^T)}+\|\nabla p\|_{L_2(\Omega^T)}\cr
&\le\Phi[T^{1/6}\|v\|_{L_2(\Omega^T)}+\|p\|_{L_2(\Omega^T)}+\varrho^*
\|v'\cdot\nabla v\|_{L_2(\Omega^T)}\cr
&\quad+\varrho^*\|v_3h\|_{L_2(\Omega^t)}+\varrho^*\|f\|_{L_2(\Omega^T)}+
\|v_0\|_{H^1(\Omega)}],\cr}
\leqno(4.26)
$$
where the third and the fourth terms we estimate by
$$\eqal{
&\|v'\cdot\nabla v\|_{L_2(\Omega^T)}\le\|v'\|_{L_{10}(\Omega^t)}
\|v\|_{W_{5/3}^{2,1}(\Omega^T)},\cr
&\|v_3h\|_{L_2(\Omega^T)}\le\|v_3\|_{W_{5/3}^{2,1}(\Omega^T)}
\|h\|_{L_{10/3}(\Omega^t)}.\cr}
\leqno(4.27)
$$
Using the energy estimate (3.10), (4.27), (4.24), (4.25) and (3.41) in (4.26) 
we obtain
$$\eqal{
&\|v\|_{W_2^{2,1}(\Omega^T)}+\|\nabla p\|_{L_2(\Omega^T)}\cr
&\le\Phi[H^2+X_1^2(T^{1/6}\|v_t\|_{L_2(\Omega^T)}^2+
\|v\|_{L_\infty(0,T;L_3(\Omega))}^2\cr
&\quad+\|\nabla v'\|_{L_2(0,T;L_3(S_2))})+G_4^2],\cr}
\leqno(4.28)
$$
where we used the estimate
$$
\|v_t\|_{L_{5/3}(\Omega^T)}\le|\Omega|^{1/6}T^{1/6}\|v_t\|_{L_2(\Omega^T)}
$$
and
$$
G_4=G_3+\|f\|_{L_2(\Omega^T)}+\|v_0\|_{H^1(\Omega)},
\leqno(4.29)
$$
where $G_3$ is defined by (4.23).

\noindent
Since $G_4\le cG$ we obtain from (4.28) inequality (4.17). This concludes 
the proof.

\noindent
Now we increase regularity from Lemma 4.3 up to 
$v\in W_2^{2+s,1+s/2}(\Omega^T)$, $s+2\ge\sigma>{5\over2}$.

\proclaim Lemma 4.4. 
Assume that $s\in\left({1\over2},1\right)$, $v$ is a weak solution to problem 
(1.1), $\varrho_*\le\varrho_0\le\varrho^*$, 
$\varrho_x(0)\in L_\infty(\Omega)$, $v\in W_2^{2+s,1+s/2}(\Omega^T)$,\break
$h\in L_\infty(0,T;L_3(\Omega))\cap L_{10\over3}(\Omega^T)$, 
$f\in W_2^{s,s/2}(\Omega^T)$, $v_0\in W_2^{1+s}(\Omega)$, 
$f\in L_2(0,T;W_{6/5}^1(\Omega))$.
Then
$$\eqal{
&\|v\|_{W_2^{2+s,1+s/2}(\Omega^T)}+\|\nabla p\|_{W_2^{s,s/2}(\Omega^T)}\cr
&\le\varphi(\varrho_*,\varrho^*,\varphi
(T^{1/2}\|v\|_{H^{2+s,1+s/2}(\Omega^T)})X_1,d_3)[\varphi
(\|v\|_{H^{2+s,1+s/2}(\Omega^T)})X_2\cr
&\quad+\|v\|_{H^{2+s,1+s/2}(\Omega^T)}X_1+\varphi(H+
\|v\|_{H^{2+s,1+s/2}(\Omega^T)}X_1+G,d_3)\cr
&\quad+K],\cr}
\leqno(4.30)
$$
where
$$\eqal{
&K=\|f\|_{W_2^{s,s/2}(\Omega^T)}+\|v_0\|_{W_2^{1+s}(\Omega)}+d_3,\cr
&X_2=X_1+\|\varrho_{xx}(0)\|_{L_q(\Omega)}+
\|\partial_t^{s/2}\varrho_x(0)\|_{L_r(\Omega)},\cr}
\leqno(4.31)
$$
where $q\le{3\over{3\over2}-s}$, $s\in(1/2,1)$, $3/2-s\le3/q\le1/2+3/r$, 
$r\le6$, $G$, $H$ are defined by (4.18).

\Proof 
From (2.6) we have
$$\eqal{
&\|v\|_{H^{s+2,s/2+1}(\Omega^T)}+\|\nabla p\|_{H^{s,s/2}(\Omega^T)}\cr
&\le\Phi[d_3+\|p\|_{H^{s,s/2}(\Omega^T)}+
\|\varrho v\cdot\nabla v\|_{H^{s,s/2}(\Omega^T)}+
\|\varrho f\|_{H^{s,s/2}(\Omega^T)}\cr
&\quad+\|v_0\|_{H^{1+s}(\Omega)}],\cr}
\leqno(4.32)
$$
where $\Phi$ is defined by (4.20).\\
Since $\varrho_x,\varrho_t\in L_\infty(\Omega^T)$ and since we are interested 
in the case $s<1$ we have
$$\eqal{
\|\varrho v\cdot\nabla v\|_{H^{s,s/2}(\Omega^T)}&\le(\varrho_x^*+\varrho_t^*)
\|v\cdot\nabla v\|_{L_2(\Omega^T)}\cr
&\quad+\varrho^*\|v\cdot\nabla v\|_{H^{s,s/2}(\Omega^T)}.\cr}
\leqno(4.33)
$$
To estimate the last term in (4.33) it is sufficient to examine the highest 
order terms. First we use the splitting
$$
\|v\cdot\nabla v\|_{H^{s,s/2}(\Omega^)}=
\|v\cdot\nabla v\|_{L_2(0,T;H^s(\Omega))}+
\|v\cdot\nabla v\|_{L_2(\Omega;H^{s/2}(0,T))}.
\leqno(4.34)
$$
It is sufficient to examine only one norm. Therefore, we consider
$$
\|v\cdot\nabla v\|_{L_2(0,T;H^s(\Omega))}=
\|\partial_x^sv\nabla v\|_{L_2(\Omega^T)}+
\|v\cdot\nabla\partial_x^sv\|_{L_2(\Omega^T)}+
\|v\cdot\nabla v\|_{L_2(\Omega^T)}.
$$
Hence we examine only the first two norms. By the H\"older inequality we have
$$\eqal{
&\|\partial_x^sv\cdot\nabla v\|_{L_2(\Omega^T)}\le
\|\partial_x^sv\|_{L_5(\Omega^T)}\|\nabla v\|_{L_{10/3}(\Omega^T)}\cr
&\le(\varepsilon_1^{1-\varkappa_1}\|v\|_{H^{2+s,1+s/2}(\Omega^T)}+
c\varepsilon_1^{-\varkappa_1}\|v\|_{L_2(\Omega^T)})
\|v\|_{H^{2,1}(\Omega^T)},\cr}
\leqno(4.35)
$$
where $\varkappa_1={{3\over2}+s\over2+s}<1$ and there is no restrictions on 
$s\in(1/2,1)$.

\noindent
Similarly, we have
$$\eqal{
&\|v\cdot\nabla\partial_x^sv\|_{L_2(\Omega^T)}\le\|v\|_{L_{10}(\Omega^T)}
\|\nabla\partial_x^sv\|_{L_{5/2}(\Omega^T)}\cr
&\le[\varepsilon_2^{1-\varkappa_1}\|v\|_{H^{2+s,1+s/2}(\Omega^T)}+
c\varepsilon_2^{-\varkappa_1}\|v\|_{L_2(\Omega^T)}]
\|v\|_{H^{2,1}(\Omega^T)}.\cr}
\leqno(4.36)
$$
Next
$$
\|\varrho f\|_{H^{s,s/2}(\Omega^T)}\le(\varrho_x^*+\varrho_t^*+\varrho^*)
\|f\|_{H^{s,s/2}(\Omega^T)}.
\leqno(4.37)
$$
Finally, we examine the term with pressure. We have
$$
\|p\|_{H^{s,s/2}(\Omega^T)}=\|p\|_{L_2(\Omega;H^{s/2}(0,T))}+
\|p\|_{L_2(0,T;H^s(\Omega))}.
\leqno(4.38)
$$
Applying $\partial_t^{s/2}$ to (3.39) and integrating the result over 
$\Omega^t$ yields
$$\eqal{
&\|\partial_t^{s/2}p\|_{L_2(\Omega^t)}\le c\bigg[\intop_0^t
(\|\partial_t^{s/2}(\nabla\varrho v_t)\|_{L_{1'}(\Omega)}^2\cr
&\quad+\|\partial_t^{s/2}(\varrho v\cdot\nabla v-\varrho f)\|_{L_{6/5}
(\Omega)}^2+\|\partial_t^{s/2}v\|_{W_{6/5}^1(\Omega)}^2)dt\bigg]^{1/2}\cr
&\le c\bigg[\intop_0^t(\|\partial_t^{s/2}\nabla\varrho\|_{L_{2'}(\Omega)}^2
\|v_t\|_{L_2(\Omega)}^2+\|\varrho_x\|_{L_{2'}(\Omega)}^2
\|\partial_t^{s/2}v_t\|_{L_2(\Omega)}^2\cr
&\quad+\|\partial_t^{s/2}\varrho v\nabla v\|_{L_{6/5}(\Omega)}^2
+\|\varrho\partial_t^{s/2}v\nabla v\|_{L_{6/5}(\Omega)}^2\cr
&\quad+\|\varrho v\partial_t^{s/2}\nabla v\|_{L_{6/5}(\Omega)}^2+
\|\partial_t^{s/2}\varrho f\|_{L_{6/5}(\Omega)}^2\cr
&\quad+\|\varrho\partial_t^{s/2}f\|_{L_{6/5}(\Omega)}^2+
\|\partial_t^{s/2}\nabla v\|_{L_{6/5}(\Omega)}^2+
\|\partial_t^{s/2}v\|_{L_{6/5}(\Omega)}^2)dt\bigg]^{1/2}\cr
&\le c\{\sup_t\|\partial_t^{s/2}\nabla\varrho\|_{L_{2'}(\Omega)}
\|v_t\|_{L_2(\Omega^t)}+\varrho_x^*\|\partial_t^{s/2}v_t\|_{L_2(\Omega^t)}\cr
&\quad+(\varrho_x^*+\varrho_t^*)[\sup_t\|v\|_{L_2(\Omega)}
\|\nabla v\|_{L_2(0,t;L_3(\Omega))}+\|f\|_{L_2(0,t;L_{6/5}(\Omega))}]\cr
&\quad+\varrho^*(\sup_t\|\partial_t^{s/2}v\|_{L_3(\Omega)}
\|\nabla v\|_{L_2(\Omega^t)}+\sup_t\|v\|_{L_2(\Omega)}
\|\partial_t^{s/2}\nabla v\|_{L_2(0,t;L_3(\Omega))}\cr
&\quad+\|\partial_t^{s/2}f\|_{L_2(0,t;L_{6/5}(\Omega))})+
\|\partial_t^{s/2}\nabla v\|_{L_2(0,t;L_{6/5}(\Omega))}\cr
&\quad+\|\partial_t^{s/2}v\|_{L_2(0,t;L_{6/5}(\Omega))}\}\equiv I_1,\cr}
$$
where $1'>1$, $2'>2$ but arbitrary close to 1 and 2, respectively.
\goodbreak

Using interpolation inequalitites (see [BIN, Ch. 3]) and the estimate for the 
weak solution we get
$$\eqal{
&\|\partial_t^{s/2}p\|_{L_2(\Omega^t)}\le I_1\le c(\sup_t
\|\partial_t^{s/2}\varrho_x\|_{L_{2'}(\Omega)}+\varrho_x^*\cr
&\quad+\varrho_t^*+\varepsilon)V_s(t)+\varphi(1/\varepsilon,d_3,\varrho^*)
(d_3+\|f\|_{L_2(0,t;L_{6/5}(\Omega))}\cr
&\quad+\|\partial_t^{s/2}f\|_{L_2(\Omega^t)}),\cr}
\leqno(4.39)
$$
where we used the notation
$$
V_s(T)=\|v\|_{H^{2+s,1+s/2}(\Omega^T)}.
$$
To estimate the norm 
$\mathop{\sup}\limits_t\|\partial_t^{s/2}\varrho_x\|_{L_r(\Omega)}$
we differentiate $(1.1)_3$ with respect to $\partial_t^{s/2}\partial_x$, 
multiply by 
$\partial_t^{s/2}\partial_x\varrho|\partial_t^{s/2}\partial_x\varrho|^{r-2}$, 
use $(1.1)_{2,4}$ and integrate over $\Omega$. Then we obtain
$$\eqal{
&{d\over dt}\|\partial_t^{s/2}\varrho_x\|_{L_r(\Omega)}\le\bigg(\intop_\Omega
|\partial_t^{s/2}v|^r|\varrho_{xx}|^rdx\bigg)^{1/r}\cr
&\quad+\varrho_x^*\|\partial_t^{s/2}v_x\|_{L_r(\Omega)}+
\|v_x\|_{L_\infty(\Omega)}\|\partial_t^{s/2}\varrho_x\|_{L_r(\Omega)}.\cr}
\leqno(4.40)
$$
Integrating (4.40) with respect to time yields
$$\eqal{
&\|\partial_t^{s/2}\varrho_x(t)\|_{L_r(\Omega)}\le\exp\bigg(\intop_0^t
\|v_x(t')\|_{L_\infty(\Omega)}dt'\bigg)\cr
&\quad\cdot\bigg[\intop_0^t\bigg(\intop_\Omega|\partial_{t'}^{s/2}v|^r
|\varrho_{xx}|^rdx\bigg)^{1/r}dt'\cr
&\quad+\varrho_x^*\intop_0^t\|\partial_t^{s/2}v_x\|_{L_r(\Omega)}dt'+
\|\partial_t^{s/2}\varrho_x(0)\|_{L_r(\Omega)}\bigg].\cr}
\leqno(4.41)
$$
On the r.h.s. of (4.41) some norm of $\varrho_{xx}$ appears. To estimate 
it we differentiate $(1.1)_3$ twice with respect to $x$, multiply by 
$\varrho_{xx}|\varrho_{xx}|^{q-2}$, use $(1.1)_{2,4}$ and integrate over 
$\Omega$. Then we obtain
$$
{d\over dt}\|\varrho_{xx}\|_{L_q(\Omega)}\le\|v_x\|_{L_\infty(\Omega)}
\|\varrho_{xx}\|_{L_q(\Omega)}+\varrho_x^*\|v_{xx}\|_{L_q(\Omega)}.
$$
Integrating the inequality with respect to time yields
$$\eqal{
&\|\varrho_{xx}(t)\|_{L_q(\Omega)}\le\exp\bigg(\intop_0^t
\|v_x(t')\|_{L_\infty(\Omega)}dt'\bigg)\cr
&\quad\cdot\bigg[\varrho_x^*\intop_0^t\|v_{xx}(t')\|_{L_q(\Omega)}dt'+
\|\varrho_{xx}(0)\|_{L_q(\Omega)}\bigg].\cr}
\leqno(4.42)
$$
Now we have to determine $r,q$ in (4.37), (4.38), respectively.

\noindent
Looking for $v\in H^{s+2,s/2+1}(\Omega^T)$ we see that 
$v_{xx}\in L_2(0,T;L_q(\Omega))$ with $q\le{3\over{3\over2}-s}$.

\noindent
Hence, (4.42) implies that $\varrho_{xx}\in L_\infty(0,T;L_q(\Omega))$,
$q\le3/\lower2pt\hbox{$(3/2-s)$}$.
Estimating the first term under the square bracket in (4.41) by
$$\eqal{
&\|\varrho_{xx}\|_{L_\infty(0,T;L_{r\lambda_1}(\Omega))}\intop_0^t
\|\partial_t^{s/2}v\|_{L_{r\lambda_2}(\Omega)}dt'\cr
&\le c(T)\|\varrho_{xx}\|_{L_\infty(0,T;L_q(\Omega))}
\|v\|_{H^{s+2,s/2+1}(\Omega^T)},\cr}
$$
we need $1/\lambda_1+1/\lambda_2=1$, $q=r\lambda_1$, 
${3\over2}-{3\over r\lambda_2}\le2$.

\noindent
Hence
$$
{3\over2}-s\le{3\over r\lambda_1}\le{1\over2}+{3\over r}
$$
which implies the restrictions
$$
1\le s+{3\over r},\quad r\le q.
\leqno(4.43)
$$

In view of the above considerations we express (4.42) in the form
$$\eqal{
&\|\varrho_{xx}(t)\|_{L_q(\Omega)}\le\exp(\|v_x\|_{L_1(0,t;L_\infty(\Omega))})
\cr
&\quad\cdot[\varrho_x^*t^{1/2}\|v\|_{H^{2+s,1+s/2}(\Omega^t)}+
\|\varrho_{xx}(0)\|_{L_q(\Omega)}],\cr}
\leqno(4.44)
$$
where
$$
q\le{3\over{3\over2}-s}.
\leqno(4.45)
$$
Using the notation
$$
V_s(T)=\|v\|_{H^{2+s,1+s/2}(\Omega^T)}
\leqno(4.46)
$$
we have
$$
\|v_x\|_{L_1(0,T;L_\infty(\Omega))}\le cT^{1/2}V_s(T),\quad {\rm for}
\ \ s>{1\over2}.
\leqno(4.47)
$$
Using (4.44) in (4.41) and using that the second term under the square 
bracket in the r.h.s. of (4.41) is estimated by
$$
c(t)V_s(t)
$$
under the assumption
$$
r\le6,
\leqno(4.48)
$$
we obtain from (4.41) the inequality
$$\eqal{
&\|\partial_t^{s/2}\varrho_x(t)\|_{L_r(\Omega)}\le\exp(t^{1/2}V_s(t))
[\exp(t^{1/2}V_s(t))\cr
&\quad\cdot(\varrho_x^*t^{1/2}V_s+\|\varrho_{xx}(0)\|_{L_q(\Omega)})t^{1/2}
V_s(t)\cr
&\quad+\varrho_x^*t^{1/2}V_s(t)+\|\partial_t^{s/2}\varrho_x(0)\|_{L_r(\Omega)}]
\cr
&\le\varphi(t^{1/2}V_s(t))[\varrho_x^*t^{1/2}V_s(t)+
\|\varrho_{xx}(0)\|_{L_q(\Omega)}+
\|\partial_t^{s/2}\varrho_x(0)\|_{L_r(\Omega)}],\cr}
\leqno(4.49)
$$
where $3/2-s\le3/q\le1/2+3/r$.

\noindent
Employing (4.49), (3.25) and (3.26) in (4.39) we obtain
$$\eqal{
&\|\partial_t^{s/2}p\|_{L_2(\Omega^t)}\le\varphi(t^{1/2}V_s(t)X_1,
\varrho_*,\varrho^*,d_3)\cr
&\quad\cdot[(X_2+\varepsilon)t^{1/2}V_s(t)+\varphi(1/\varepsilon,d_3,
\varrho_*,\varrho^*)(\|f\|_{L_2(0,t;L_{6/5}(\Omega))}\cr
&\quad+\|\partial_t^{s/2}f\|_{L_2(\Omega^t)})].\cr}
\leqno(4.50)
$$
Now we examine $\|p\|_{L_2(0,t;H^s(\Omega))}$. Applying $\partial_x^s$ to 
(3.39) and taking the $L_2(\Omega)$ norm we obtain
$$\eqal{
&\|\partial_x^sp\|_{L_2(\Omega)}\le c(\|p_x\|_{L_2(\Omega)}+
\|p\|_{L_2(\Omega)})\le c(\|\nabla\varrho v_t\|_{L_{6/5}(\Omega)}\cr
&\quad+\|\varrho v\nabla v\|_{L_2(\Omega)}+\|\varrho f\|_{L_2(\Omega)}+
\|v\|_{W_2^{2-1/2}(S)})\cr
&\le c\varrho_x^*\|v_t\|_{L_{6/5}(\Omega)}+c\varrho^*\|v\|_{L_3(\Omega)}
\|\nabla v\|_{L_2(\Omega)}\cr
&\quad+c\varrho^*\|f\|_{L_2(\Omega)}+c\|v\|_{W_2^2(\Omega)}.\cr}
\leqno(4.51)
$$
We obtain from (4.51) after integration with respect to time the inequality
$$\eqal{
\|p_{,x}\|_{L_2(\Omega^t)}&\le c\varrho_x^*V_s(t)+\varphi(d_3,\varrho_x^*)
(\varepsilon V_s(t)\cr
&\quad+c(1/\varepsilon)d_3+\|f\|_{L_2(\Omega^t)}).\cr}
\leqno(4.52)
$$
Using (4.33)--(4.37), (4.50) and (4.52) in (4.32) implies the inequality
$$\eqal{
&\|v\|_{H^{2+s,1+s/2}(\Omega^T)}+\|\nabla p\|_{H^{s,s/2}(\Omega^T)}\cr
&\le\varphi(T^{1/2}V_s(T)X_1,\varrho_*,\varrho^*,d_3)
[(X_2+\varepsilon)T^{1/2}V_s(T)\cr
&\quad+\varphi(1/\varepsilon,d_3)(\|v\|_{H^{2,1}(\Omega^T)}+
\|f\|_{H^{s,s/2}(\Omega^T)})+\varphi(1/\varepsilon,\varrho_*,\varrho^*,d_3)].
\cr}
\leqno(4.53)
$$
For sufficiently small $\varepsilon$ and (4.17) we obtain (4.30). This 
concludes the proof.

\Remark{4.5.} 
In formulas (3.25), (3.26), (4.41) and (4.44) we have the expression
$$
I(t)=\exp\intop_0^t\|v_x(t')\|_{L_\infty(\Omega)}dt'.
$$
In view of the imbedding
$$
\|v\|_{L_2(0,T;W_\infty^1(\Omega))}\le cV_s(T),\quad s>{1\over2}
$$
we obtain the estimate
$$
I(t)\le\exp(ct^{1/2}V_s(t)),\quad t\le T,
$$
which is not convenient because factor $t^{1/2}$ appears under the exponent 
functions. The difficulty can be cancelled in virtue of the assumption
$$
1\le\|v\|_{W_\infty^1(\Omega)}.
$$
The above inequality is not restrictive because the case 
$\|v\|_{W_\infty^1(\Omega)}\le1$ implies in view of Lemma 3.5 the following 
estimate for solutions to problem~(1.1)
$$\eqal{
&\|v\|_{W_2^{2,1}(\Omega^T)}+\|\nabla p\|_{L_2(\Omega^T)}\le\varphi
(\varrho_*,\varrho^*,d_3,e^{cT})\cr
&\quad\cdot[\|f\|_{L_2(\Omega^T)}+\|v(0)\|_{H^1(\Omega)}],\quad c>0.\cr}
$$
Hence regularity $H^{s+2,s/2+1}(\Omega^T)$ follows immediately. In this 
case the constant in the above estimate depends on time but this does not 
imply any restrictions on magnitudes of the data.

\noindent
Now we pass to problem (3.17). Then we have

\proclaim Lemma 4.6. 
Assume that ${1\over2}<\sigma\le s<1$ and $\sigma$ can be chosen as very 
close to $s$.\\
Let us take Remark 4.5 under account and let
$$\eqal{
B&=\|g\|_{L_2(0,T;L_{6/5}(\Omega))}+\|f_3\|_{L_2(0,T;L_{4/3}(S_2))}\cr
&\quad+\|h(0)\|_{L_2(\Omega)}<\infty,\cr}
\leqno(4.54)
$$
be sufficiently small.\\
Let $f,g\in L_2(\Omega^T)$, $h(0)\in H^1(\Omega)$, 
$v\in H^{2+\sigma,1+\sigma/2}(\Omega^T)$ and let $(v,\varrho)$ be the weak 
solution to (1.1) described by Lemma 3.5.
Then solutions to (3.17) satisfy the inequality
$$\eqal{
&\|h\|_{H^{2,1}(\Omega^T)}+\|\nabla q\|_{L_2(\Omega^T)}\le\varphi
(\varrho_*,\varrho^*,\varphi(V_\sigma(T))X_1,d_3)\cr
&\quad\cdot[\varphi(V_\sigma(T))(X_1+B)+K_0(T)],\cr}
\leqno(4.55)
$$
where $X_1$ is introduced in (3.22) and
$$
K_0(T)=\|g\|_{L_2(\Omega^T)}+\|f\|_{L_2(\Omega^T)}+\|h(0)\|_{H^1(\Omega)}.
\leqno(4.56)
$$

\Proof 
For solutions to (3.17) we have
$$\eqal{
&\|h\|_{H^{2,1}(\Omega^T)}+\|\nabla q\|_{L_2(\Omega^T)}\le\varphi
(\varrho_*,\varrho^*,\varphi(V_\sigma(T))X_1,d_3)\cr
&\quad\cdot[\|v\cdot\nabla h\|_{L_2(\Omega^T)}+
\|h\cdot\nabla v\|_{L_2(\Omega^T)}+\|g\|_{L_2(\Omega^T)}\cr
&\quad+\varphi(V_\sigma(T)X_1)X_1(\|v_t\|_{L_2(\Omega^T)}+
\|v\cdot\nabla v\|_{L_2(\Omega^T)}+\|f\|_{L_2(\Omega^T)})\cr
&\quad+\|h(0)\|_{H^1(\Omega)}].\cr}
\leqno(4.57)
$$
We need the inequalities
$$\eqal{
\|v\cdot\nabla h\|_{L_2(\Omega^T)}&\le\|\nabla h\|_{L_3(\Omega^T)}
\|v\|_{L_6(\Omega^T)}\cr
&\le\varepsilon\|h\|_{H^{2,1}(\Omega^T)}+\varphi(\|v\|_{W_2^{2,1}(\Omega^T)})
\|h\|_{L_2(\Omega^T)},\cr}
\leqno(4.58)
$$
$$\eqal{
\|h\cdot\nabla v\|_{L_2(\Omega^T)}&\le\|h\|_{L_6(\Omega^T)}
\|\nabla v\|_{L_3(\Omega^T)}\cr
&\le\varepsilon\|h\|_{H^{2,1}(\Omega^T)}+\varphi(\|v\|_{W_2^{2,1}(\Omega^T)})
\|h\|_{L_2(\Omega^T)}\cr}
\leqno(4.59)
$$
$$
\|v\cdot\nabla v\|_{L_2(\Omega^T)}\le\|v\|_{L_\infty(\Omega^T)}
\|\nabla v\|_{L_2(\Omega^T)}\le d_3\|v\|_{H^{2+\sigma,1+\sigma/2}(\Omega^T)}.
\leqno(4.60)
$$
In view of (4.6) we have
$$\eqal{
\|h\|_{L_2(\Omega^T)}&\le\varphi(\varrho_*,\varrho^*,d_3)\exp(V_\sigma(T))[B\cr
&\quad+X_1(V_\sigma(T)+\|f\|_{L_2(0,T;L_{4/3}(S_2))})].\cr}
\leqno(4.61)
$$
Employing (4.58)--(4.61) in (4.57) implies (4.55).This concludes the proof.

\noindent
Finally, we have

\proclaim Theorem 4.7. 
Let assumptions of Lemmas 4.4, 4.6 hold. Let us take Remark 4.5 under 
account. Let
$$
X=X_1+X_2+B
\leqno(4.68)
$$
Then for sufficiently small $X$ the estimate holds
$$
V_s(T)\le\varphi(\varrho_*,\varrho^*,G,K,K_0,X,L),\quad {1\over2}<s<1,
\leqno(4.69)
$$
where $K_0$ is introduced in (4.56), $K$ in(4.31), $G$ in (4.18) and
$$
L=\|f_3\|_{L_2(0,T;H^{1/2}(S))}+\|f\|_{H^{s,s/2}(\Omega^T)}.
\leqno(4.70)
$$

\Proof Applying a fixed point argument in (4.30) we obtain 
for sufficiently small $X_2$ the inequality
$$
V_s(T)\le\varphi(\varrho_*,\varrho^*,\|h\|_{H^{2,1}(\Omega^T)},G,K,X_2).
\leqno(4.71)
$$
Using (4.55) with $\sigma=s$ in (4.71) and applying again a fixed point 
argument for sufficiently small $X$ we obtain (4.69). This concludes the proof.

\section{5. Existence}

We prove the existence of solutions to problem (1.1) by the Leray\--Schauder 
fixed point theorem. For this purpose we construct a mapping $\Phi$ in the 
following way.

\noindent
Let $\tilde v\in W_2^{2+\sigma,1+\sigma/2}(\Omega^T)$, $\divv\tilde v=0$, 
$\tilde v\cdot\bar n|_S=0$, $\sigma\in(1/2,1)$ be given.

\noindent
Then $\varrho=\varrho(\tilde v)$ is a solution to the problem
$$
\varrho_t+\tilde v\cdot\nabla\varrho=0,\quad \varrho|_{t=0}=\varrho_0.
\leqno(5.1)
$$
Let
$$
v=\Phi(\tilde v,\lambda)
\leqno(5.2)
$$
be a solution to the problem
$$\eqal{
&[\varrho_0(1-\lambda)+\lambda\varrho(\tilde v)]v_t-\divv\T(v,p)\cr
&=\lambda[-\varrho(\tilde v)\tilde v\cdot\nabla\tilde v+\varrho(\tilde v)f]\cr
&\divv v=0\cr
&v\cdot\bar n|_S=0,\ \ 
\bar n\cdot\D(v)\cdot\bar\tau_\alpha+\gamma v\cdot\bar\tau_\alpha|_S=0,\ \ 
\alpha=1,2,\cr
&v|_{t=0}=v_0,\cr}
\leqno(5.3)
$$
where $\lambda\in[0,1]$.

In view of Lemma 4.4, Remark 4.5, Lemma 4.6 and Theorem 4.7 we have

\proclaim Lemma 5.1. 
Assume that $v_0\in H^{1+s}(\Omega)$, $v_{0,x_3}\in H^1(\Omega)$, 
$f\in H^{s,s/2}(\Omega^T)\break\cap L_2(0,T;W_{6/5}^1(\Omega))$, 
$f_{,x_3}\in L_2(\Omega^T)$, $s\in(1/2,1)$, $\varrho_0\in W_q^2(\Omega)$, 
$3<q\le{3\over3/2-s}$. Assume also that there exist positive constants 
$\varrho_*<\varrho^*$ such that $\varrho_*\le\varrho_0\le\varrho^*$.
Then the mapping (5.2) has a fixed point belonging to 
$H^{2+s,1+s/2}(\Omega^T)$.

\Remark{5.2.}
In Lemma 4.4 there is assumption that 
$\partial_t^{s/2}\varrho_x|_{t=0}\in L_r(\Omega)$, where $r\le6$ and is such 
that $1\le s+{3\over r}$ for $s\in(1/2,1)$, $r\le q$.

\noindent
We see that the condition is satisfied in view of the relations
$$
\varrho_0\in W_q^2(\Omega),\quad r\le q,\quad v_0\in H^{1+s}(\Omega).
$$
We have
$$\eqal{
&\|\partial_t^{s/2}\varrho_x|_{t=0}\|_{L_r(\Omega)}\le c
(\|\partial_t\varrho_x|_{t=0}\|_{L_r(\Omega)}+
\|\varrho_x|_{t=0}\|_{L_r(\Omega)})\cr
&\le c\|(v_x\cdot\nabla\varrho+v\cdot\nabla\varrho_x)|_{t=0}\|_{L_r(\Omega)}+
c\|\varrho_{0,x}\|_{L_r(\Omega)}\cr
&\le c(\|v_0\|_{H^{1+s}(\Omega)}\|\varrho_{0,x}\|_{W_q^1(\Omega)}+
\|\varrho_{0,x}\|_{L_q(\Omega)}).\cr}
$$

\noindent
{\bf Proof of Lemma 5.1.} 
In view of the a priori estimate (4.69) we have toexamine the other 
assumptions of the Leray-Schauder fixed point theorem.

\noindent
For $\lambda=0$ we have the existence of a unique solution.

\noindent
We have that
$$
\varrho_*\le\varrho_0(1-\lambda)+\lambda\varrho(\tilde v)\le\varrho^*
$$
since $\varrho_0$ is continuous and $\varrho_0\in W_p^1(\Omega)$, $p>3$ 
we have that $\varrho=\varrho(\tilde v)$ is continuous because
$$\eqal{
&\|\varrho_x(t)\|_{L_p(\Omega)}\le\exp
(\|\tilde v\|_{L_1(0,t;L_\infty(\Omega))})\|\varrho_x(0)\|_{L_p(\Omega)},\cr
&\|\varrho_t(t)\|_{L_p(\Omega)}\le\|\tilde v\|_{L_\infty(\Omega^t)}\sup_t
\|\varrho_x\|_{L_p(\Omega)}\cr
&\le\|\tilde v\|_{L_\infty(\Omega^t)}\exp
(\|\tilde v\|_{L_1(0,t;L_\infty(\Omega))})\|\varrho_x(0)\|_{L_p(\Omega)}.\cr}
\leqno(5.4)
$$
Then
$$\eqal{
&\sup_{x,x',t}{|\varrho(x,t)-\varrho(x',t)|\over|x-x'|^{1/p'}}\le\sup
{|\varrho(x_1,x_2,x_3,t)-\varrho(x'_1,x_2,x_3,t)|\over|x_1-x'_1|^{1/p'}}\cr
&\quad+\sup{|\varrho(x'_1,x_2,x_3,t)-\varrho(x'_1,x'_2,x_3,t)|\over
|x_2-x'_2|^{1/p'}}\cr
&\quad+\sup
{|\varrho(x'_1,x'_2,x_3,t)-\varrho(x'_1,x'_2,x'_3,t)|\over|x_3-x'_3|^{1/p'}}\cr
&\le\sup\bigg(\intop_{x'_1}^{x_1}|\varrho_z(z,x_2,x_3,t)|^pdz\bigg)^{1/p}+
\sup\bigg(\intop_{x'_2}^{x_2}|\varrho_z(x_1,z,x_3,t)|^pdz\bigg)^{1/p}\cr
&\quad+\sup\bigg(\intop_{x'_3}^{x_3}|\varrho_z(x_1,x_2,z,t)|^pdz\bigg)^{1/p}
\le c\sup_t\|\varrho_x\|_{L_\infty(\Omega)},\cr}
\leqno(5.5)
$$
where $x_i'<x_i$, $i=1,2,3.$ Moreover,
$$
{|\varrho(x,t)-\varrho(x,t')|\over|t-t'|^{1/p'}}\le
\bigg(\intop_{t'}^t|\varrho_t|^pdt\bigg)^{1/p}.
\leqno(5.6)
$$
Hence we have continuity of $\varrho$ if $\varrho_x(0)\in L_\infty(\Omega)$.

\noindent
Assuming that $\tilde v\in W_2^{2+\sigma,1+\sigma/2}(\Omega^T)$, 
$\sigma\in(1/2,1)$ we show that
$$
\|\varrho\tilde v\cdot\nabla\tilde v\|_{W_2^{s,s/2}(\Omega^T)}\le\varphi
(\|\tilde v\|_{W_2^{2+\sigma,1+\sigma/2}(\Omega^T)}),
$$
where $1>s>\sigma$.

\noindent
We have to assume that $\sigma>1/2$ because we need the imbeddings
$$
\|\tilde v\|_{L_\infty(\Omega^T)}+\|\tilde v_x\|_{L_2(0,T;L_\infty(\Omega))}
\le c\|\tilde v\|_{W_2^{2+\sigma,1+\sigma/2}(\Omega^T)}.
$$
Hence, we have shown that
$$
\Phi:\ W_2^{2+\sigma,1+\sigma/2}(\Omega^T)\times[0,1]\to
W_2^{2+s,1+s/2}(\Omega^T)\subset W_2^{2+\sigma,1+\sigma/2}(\Omega^T)
\leqno(5.7)
$$
where the last imbedding is compact. Then mapping $\Phi$ is compact.

\noindent
To show continuity of mapping $\Phi$ we introduce the differences
$$
\tilde V=\tilde v_1-\tilde v_2,\quad V=v_1-v_2,\quad P=p_1-p_2.
\leqno(5.8)
$$
Then $(V,P)$ is a solution to the problem
$$\eqal{
&[\varrho_0(1-\lambda)+\lambda\varrho(\tilde v_2)]V_t-\divv\T(V,P)\cr
&=-\lambda(\varrho(\tilde v_1)-\varrho(\tilde v_2))v_{1t}-
(\varrho(\tilde v_1)\tilde v_1\cdot\nabla\tilde v_1-\varrho(\tilde v_2)
\tilde v_2\cdot\nabla\tilde v_2)\cr
&\quad+(\varrho(\tilde v_1)-\varrho(\tilde v_2))f,\cr
&\divv V=0,\cr
&V\cdot\bar n|_S=0,\ \ 
\bar n\cdot\D(V)\cdot\bar\tau_\alpha+\delta_{1j}\gamma V\cdot
\bar\tau_\alpha|_{S_j}=0,\ \ \alpha=1,2,\ \ j=1,2,\cr
&V|_{t=0}=0.\cr}
\leqno(5.9)
$$
From (5.9) we have
$$
\|V\|_{W_2^{2+s,1+s/2}(\Omega^T)}\le\varphi(A)
\|\tilde V\|_{W_2^{2+\sigma,1+\sigma/2}(\Omega^T)},
\leqno(5.10)
$$
where
$$
\sum_{i=1}^2\|\tilde v_i\|_{W_2^{2+s,1+s/2}(\Omega^T)}\le A,
$$
where the last estimate is shown in Section 4.

\noindent
Let us examine continuity of $\Phi$ with respect to $\lambda$. For this 
purpose we examine
$$\eqal{
&[\varrho_0(1-\lambda_i)+\lambda_i\varrho(\tilde v)]v_{it}-\divv\T(v_i,p_i)\cr
&=\lambda_i[-\varrho(\tilde v)\tilde v\cdot\nabla\tilde v+\varrho(\tilde v)f],
\cr
&\divv v_i=0\cr
&\bar n\cdot v_i|_S=0,\ \ 
\bar n\cdot\D(v_i)\cdot\bar\tau_\alpha+\delta_{1j}\gamma v_i\cdot
\bar\tau_\alpha|_{S_j}=0,\ \ \alpha=1,2,\ \ j=1,2,\cr
&v_i|_{t=0}=v_0,\cr}
\leqno(5.11)
$$
where $i=1,2$.

\noindent
Introducing the differences
$$
V=v_1-v_2,\quad P=p_1-p_2,\quad \Lambda=\lambda_1-\lambda_2
\leqno(5.12)
$$
we see that they satisfy the problem
$$\eqal{
&[\varrho_0(1-\lambda_2)+\lambda_2\varrho(\tilde v)]V_t-\divv\T(V,P)\cr
&=\Lambda(\varrho_0-\varrho(\tilde v))v_{1t}+\Lambda(-\varrho(\tilde v)
\tilde v\cdot\nabla\tilde v+\varrho(\tilde v)f),\cr
&\divv =0\cr
&V\cdot\bar n|_S=0,\ \ 
\nu\bar n\cdot\D(V)\cdot\bar\tau_\alpha+\delta_{1j}\gamma V\cdot
\bar\tau_\alpha|_{S_j}=0,\ \ \alpha=1,2,\ \ j=1,2,\cr
&V|_{t=0}=0.\cr}
\leqno(5.13)
$$
For solutions to (5.13) we have
$$\eqal{
&\|V\|_{W_2^{2+s,1+s/2}(\Omega^T)}+\|\nabla P\|_{W_2^{s,s/2}(\Omega^T)}\cr
&\le\varphi(A)(1+\|v_{1t}\|_{W_2^{s,s/2}(\Omega^T)})\Lambda.\cr}
\leqno(5.14)
$$
Hence, the continuity with respect to $\lambda$ follows.

\noindent
Applying the Leray-Schauder fixed point theorem we prove Lemma~5.1. This 
concludes the proof.

Now we prove uniqueness

\proclaim Lemma 5.3. 
Assume that $\varrho_*\le\varrho\le\varrho^*$, 
$\varrho\in L_2(0,T;H^2(\Omega))$, $v\in H^{2,1}(\Omega^T)$, 
$f\in L_2(\Omega^T)$.
Then we have uniqueness of solutions to problem (1.1).

\Proof 
Assume that we have two solutions $(v_i,\varrho_i,p_i)$, $i=1,2$, to 
problem (1.1). Let
$$
V=v_1-v_2,\quad R=\varrho_1-\varrho_2,\quad P=p_1-p_2.
\leqno(5.15)
$$
Then functions (5.15) are solutions to the problem
$$\eqal{
&\varrho_1V_t+\varrho_1v_1\cdot\nabla V+\varrho_1V\cdot\nabla v_2-\divv\T(V,P)
\cr
&=Rf-R(v_{2t}+v_2\cdot\nabla v_2),\cr
&\divv V=0,\cr
&\varrho_1R_t+\varrho_1v_1\cdot\nabla R+\varrho_1V\cdot\nabla\varrho_2=0,\cr
&V\cdot\bar n|_S=0,\ \ 
(\nu\bar n\cdot\T(V,P)\cdot\bar\tau_\alpha+\delta_{1j}\gamma V\cdot
\bar\tau_\alpha)|_{S_j}=0,\ \ \alpha=1,2,\ \ j=1,2,\cr
&V|_{t=0}=0,\ \ R|_{t=0}=0.\cr}
\leqno(5.16)
$$
Multiplying $(5.16)_1$ by $V$, integrating over $\Omega$, and employing the 
equation of continuit for $\varrho_1$ and using the Korn inequality we obtain
$$\eqal{
&{1\over2}{d\over dt}\intop_\Omega\varrho_1V^2dx+\|V\|_{H^1(\Omega)}^2+
\gamma\intop_{S_1}|V\cdot\bar\tau_\alpha|^2dS_1\cr
&\le-\intop_\Omega V\cdot\nabla v_2\cdot Vdx+\intop_\Omega Rf\cdot Vdx-
\intop_\Omega R(v_{2t}+v_2\cdot\nabla v_2)\cdot Vdx.\cr}
\leqno(5.17)
$$
The first term on the r.h.s. of (5.17) is estimated by
$$
\varepsilon\|V\|_{L_6(\Omega)}^2+c(1/\varepsilon)
\|\nabla v_2\|_{L_3(\Omega)}^2\|V\|_{L_2(\Omega)}^2,
$$
the second by
$$
\varepsilon\|V\|_{L_6(\Omega)}^2+c(1/\varepsilon)\|f\|_{L_2(\Omega)}^2
\|R\|_{L_3(\Omega)}^2
$$
and the last by
$$
\varepsilon\|V\|_{L_6(\Omega)}^2+c(1/\varepsilon)(\|v_{2t}\|_{L_2(\Omega)}^2+
\|v_2\|_{L_4(\Omega)}^2\|\nabla v_2\|_{L_4(\Omega)}^2)\|R\|_{L_3(\Omega)}^2.
$$
Using the estimates in (5.17) and assuming that $\varepsilon$ is sufficiently 
small we arrive to the inequality
$$\eqal{
&{d\over dt}\intop_\Omega\varrho_1V^2dx+\|V\|_{H^1(\Omega)}^2\le c
\|\nabla v_2\|_{L_3(\Omega)}^2\|V\|_{L_2(\Omega)}^2\cr
&\quad+c(\|f\|_{L_2(\Omega)}^2+\|v_{2t}\|_{L_2(\Omega)}^2+
\|v_2\|_{L_4(\Omega)}^2\|\nabla v_2\|_{L_4(\Omega)}^2)\|R\|_{L_3(\Omega)}^2.
\cr}
\leqno(5.18)
$$
Multiplying $(5.16)_3$ by $R^2$, integrating over $\Omega$ and using the 
equation of continuity for $\varrho_1$ yields
$$
{1\over3}{d\over dt}\intop_\Omega\varrho_1R^3dx=-\intop_\Omega\varrho_1V\cdot
\nabla\varrho_2R^2dx.
\leqno(5.19)
$$
Estimating the r.h.s. by
$$\eqal{
&\bigg(\intop_\Omega|V\cdot\nabla\varrho_2|^3dx\bigg)^{1/3}
\bigg(\intop_\Omega(\varrho_1R^2)^{3/2}dx\bigg)^{2/3}\cr
&\le\varphi(\varrho^*)\|V\|_{L_6(\Omega)}\|\nabla\varrho_2\|_{L_6(\Omega)}
\bigg(\intop_\Omega\varrho_1R^3dx\bigg)^{2/3}\cr}
$$
we obtain from (5.19) the inequality
$$\eqal{
&{d\over dt}\bigg(\intop_\Omega\varrho_1R^3dx\bigg)^{2/3}\le\varphi(\varrho^*)
\|V\|_{L_6(\Omega)}\|\nabla\varrho_2\|_{L_6(\Omega)}
\bigg(\intop_\Omega\varrho_1R^3dx\bigg)^{1/3}\cr
&\le\varepsilon\|V\|_{L_6(\Omega)}^2+\varphi(1/\varepsilon,\varrho^*)
\|\nabla\varrho_2\|_{L_6(\Omega)}^2
\bigg(\intop_\Omega\varrho_1R^3dx\bigg)^{2/3}\cr}
\leqno(5.20)
$$
Adding (5.18) and (5.20) with $\varepsilon$ sufficiently small and defining 
the quantities
$$\eqal{
&Y(t)=\intop_\Omega\varrho_1V^2dx+
\bigg(\intop_\Omega\varrho_1R^3dx\bigg)^{2/3},\cr
&A(t)=\|\nabla\varrho_2\|_{L_6(\Omega)}^2+\|v_{2t}\|_{L_2(\Omega)}^2+
\|v_2\|_{L_4(\Omega)}^2\|\nabla v_2\|_{L_4(\Omega)}^2+
\|f\|_{L_2(\Omega)}^2\cr}
$$
we obtain
$$
{d\over dt}Y\le\varphi(\varrho^*,\varrho_*)AY
\leqno(5.21)
$$
Hence for $I=\intop_0^TAdt<\infty$ we have uniqueness.

\noindent
Moreover, we see that
$$\eqal{
I&\le\|\varrho_2\|_{L_2(0,T;H^2(\Omega))}^2+
\|v_2\|_{H^{2,1}(\Omega^T)}^2(1+\|v_2\|_{H^{2,1}(\Omega^T)}^2)\cr
&\quad+\|f\|_{L_2(\Omega^T)}<\infty.\cr}
$$
This concludes the proof.

\Remark{5.4.}
The assumptions of Lemma 5.3 are satisfied in view of the assumptions of 
Lemma 5.1.

\section{References}

\item{[Z1]} Zaj\c aczkowski, W. M.: Global existence of axially symmetric 
solutions of incompressible Navier-Stokes equations with large angular 
component of velocity, Colloq. Math. 100 (2004), 243--263.

\item{[RZ]} Renc\l awowicz, J.; Zaj\c aczkowski, W. M.: Large time regular 
solutions to the Navier-Stokes equations in cylindrical domains, Top. Meth. 
Nonlin. Anal. 32 (2008), 69--87.

\item{[Z2]} Zaj\c aczkowski, W. M.: Long time existence of regular solutions 
to the Navier-Stokes equations in cylindrical domains under boundary slip 
conditions, Studia Math. 169 (3) (2005), 243--285.

\item{[BIN]} Besov, O. V.; Il'in, V. P.; Nikolskii, S. M.: Integral 
representation of functions and theorems of imbedding, Nauka, Moscow 1975 
(in Russian).

\item{[BZ]} Burnat, M.; Zaj\c aczkowski, W. M.: On local motion of 
a~compressible barotropic viscous fluid with the boundary slip condition, 
Top. Meth. Nonlin. Anal. 10 (1997), 195--223.

\item{[Z3]} Zaj\c aczkowski, W. M.: Nonstationary Stokes system in 
Sobolev-\break-Slobodetski spaces

\item{[Z4]} Zaj\c aczkowski, W. M.: On global regular solutions to the 
Navier-Stokes equations in cylindrical domains Top. Meth. Nonlin. Anal. 
(2011).

\item{[AKM]} Antontzev, S. N.; Kazhikhov, A. V.; Monakhov, V. N.: 
Boundary problems for mechanics of nonhomogeneous fluids, Nauka, Novosibirsk 
1983 (in Russian).

\item{[NZ1]} Nowakowski, B.; Zaj\c aczkowski, W. M.: Global existence of 
soltuins to Navier-Stokes equations in cylindrical domains, Appl. Math. 
36 (2) (2009) 169--182.

\item{[NZ2]} Nowakowski, B.; Zaj\c aczkowski, W. M.: Global attractor for 
Navier\--Stokes equaitons in cylindrical domains, Appl. Math. 36 (2) (2009) 
183--194.

\item{[Z5]} Zaj\c aczkowski, W. M.: Global special regular solutions to the 
Navier\--Stokes equations in a cylindrical domain without the axis of 
symmetry, Top. Meth. Nonlin. Anal. 24 (2004) 69--105.

\item{[Z6]} Zaj\c aczkowski, W. M.: Global regular solutions to the 
Navier-Stokes equations in a cylinder, Banach Center Publ. 74 (2006), 235--255.

\item{[A1]} Alame, W.: On the existence of solutions for the nonstationary\break
Stokes system with slip boundary conditions, Appl. Math. 32 (2005), 195--223.

\item{[P]} Lions, P. L.: Mathematical topics in fluid mechanics, V.1 -- 
Incompressible models, Clarendon Press-Oxford, 1996.

\bye